\documentclass{article}
\usepackage{amssymb}
\usepackage{graphicx}
\newtheorem{theorem}{Theorem}
\newtheorem{prop}{Proposition}
\newtheorem{lemma}{Lemma}
\newtheorem{sublemma}{Sublemma}
\newtheorem{corr}{Corollary}

\begin{document}

\title{Decay of correlations for non-H\"older observables}
\author{Vincent Lynch}
\maketitle

\begin{abstract}
We consider the general question of estimating decay of correlations
for non-uniformly expanding maps, for classes of observables which
are much larger than the usual class of H\"older continuous functions.
Our results give new estimates for many non-uniformly expanding
systems, including Manneville-Pomeau maps, many one-dimensional
systems with critical points, and \emph{Viana maps}.  In many
situations, we also obtain a Central Limit Theorem for a much larger
class of observables than usual.

Our main tool is an extension of the coupling method introduced by
L.-S. Young for estimating rates of mixing on certain non-uniformly
expanding \emph{tower} maps.
\end{abstract}

\section{Introduction}

In this paper, we are interested in mixing properties (in particular,
\emph{decay of correlations}) of non-uniformly expanding maps.  Much
progress has been made in recent years, with upper estimates being obtained
for many examples of such systems.  Almost invariably, these estimates are
for observables which are H\"older continuous.  Our aim here is to extend
the study to much larger classes of observables.

Let $f:(X,\nu) \circlearrowleft$ be some mixing system.
We define a \emph{correlation function}

\[ \mathcal{C}_n (\varphi,\psi;\nu) = \left| \int (\varphi \circ f^n)
   \psi d\nu - \int \varphi d\nu \int \psi d\nu \right| \]
for $\varphi,\psi \in \mathcal{L}^2$.  The rate at which this sequence
decays to zero is a measure of how quickly $\varphi \circ f^n$ becomes
independent from $\psi$.  It is well known that for any non-trivial mixing
system, there exist $\varphi,\psi \in \mathcal{L}^2$ for which correlations
decay arbitrarily slowly.  For this reason, we must restrict at least one
of the observables to some smaller class of functions, in order to get an
upper bound for $\mathcal{C}_n$.

Here, we present a result which is general in the context of \emph{towers},
as introduced by L.-S. Young (\cite{Y2}).  There are many examples of
systems which admit such towers, and we shall see that under a fairly weak
assumption on the relationship between the tower and the system (which is
satisfied in all the examples we mention) we get estimates for certain
classes of observables with respect to the system itself.  One of the main
strengths of this method is that these classes of observables may be defined
purely in terms of their regularity with respect to the manifold; this
contrasts with some results, where regularity is considered with respect to
some Markov partition.

\paragraph{Acknowledgements.}  This paper is based on work contained in my
Ph.D. thesis, written at the University of Warwick.  Many thanks to my
supervisor, Stefano Luzzatto, and also to Omri Sarig and Peter Walters.
I am grateful for many conversations with Mark Holland, Frederic Paccaut,
Gavin Band and Mike Todd.

\section{Statement of results}

\label{obsintro}

Let us start by defining the classes of observables we consider.
Let $(X,d)$ be a metric space.
For a given function $\psi:X \rightarrow \mathbb{R}$, for each
$\varepsilon>0$ we write
\[ \mathcal{R}_\varepsilon (\psi) := \sup \{ |\psi(x) - \psi(y)| :
   d(x,y) \leq \varepsilon \} .\]
We see that $\mathcal{R}_\varepsilon (\psi) \rightarrow 0$ as $\varepsilon
\rightarrow 0$ if and only if $\psi$ is uniformly continuous.  We define
classes of functions corresponding to different rates of decay for
$\mathcal{R}_\varepsilon (\psi)$, as follows:

\emph{Class $(R1,\gamma)$, $\gamma \in (0,1]$:} $\psi \in (R1,\gamma)$ if
$\mathcal{R}_\varepsilon (\psi) = \mathcal{O}(\varepsilon^\gamma)$.

\emph{Class $(R2,\gamma)$, $\gamma \in (0,1)$:} $\psi \in (R2,\gamma)$ if
$\mathcal{R}_\varepsilon (\psi) = \mathcal{O}( \exp \{ - | \log \varepsilon
|^\gamma \})$.

\emph{Class $(R3,\gamma)$, $\gamma>1$:} $\psi \in (R3,\gamma)$ if
$\mathcal{R}_\varepsilon (\psi) = \mathcal{O}( \exp \{ - (\log | \log
\varepsilon|)^\gamma \})$.

\emph{Class $(R4,\gamma)$, $\gamma>1$:} $\psi \in (R4,\gamma)$ if
$\mathcal{R}_\varepsilon (\psi) = \mathcal{O}( | \log \varepsilon
|^{-\gamma} )$.

We write $(Rn) = \cup_\gamma (Rn,\gamma)$ for each $n$.  Note $(R1)$
is the class of functions which are H\"older continuous.  Also note
that $(R1) \subset (R2) \subset (R3) \subset (R4)$ with the inclusions
strict, and in fact $(Rn,\gamma) \supset (R(n+1))$ for every
$n=1,2,3$, with $\gamma$ in the appropriate interval.  It is not difficult
to find observables with exactly any of these regularities; for instance,
the function $\psi:[0,\frac{1}{2}] \rightarrow \mathbb{R}$ given by
\begin{eqnarray*} \psi(x) = \left\{ \begin{array}{rcl} | \log x|^{-\gamma}
  && x>0  \\  0 && x=0 \end{array} \right. \end{eqnarray*}
for some $\gamma>1$, has the property that $\mathcal{R}_\varepsilon (\psi)
= \psi (\varepsilon)$ for all small $\varepsilon$.

Before we state the exact technical result, we shall illustrate its
implications by stating results for a number of different classes of system.
In each case, the estimate given for H\"older continuous observables (i.e.
case (R1)) has been obtained previously, but we are able to give
estimates for observables in each of the other classes, most of which are
new.  Often the same estimates given in the H\"older case also apply to one
of the larger classes; it seems that when the rate of mixing for H\"older
observables is sub-exponential, we can expect the same rate of mixing to hold
for a larger class of observables.  For instance, when H\"older observables
give polynomial mixing, we tend to get the same speed of mixing for
observables in some class (R4,$\gamma$), $\gamma>1$.

All of our results shall take the following form.  Given a system
$f:X \circlearrowleft$, a mixing acip $\nu$, and $\varphi \in
\mathcal{L}^\infty (X,\nu)$, $\psi \in \mathcal{I}$, for some class
$\mathcal{I}=(Ri,\gamma)$ as above, we obtain in each example an
estimate of the form
\[ \mathcal{C}_n (\varphi,\psi;\nu) \leq \| \varphi \|_\infty C(\psi)
   u_n ,\]
where $\| \cdot \|_\infty$ is the usual norm on $\mathcal{L}^\infty(X,\nu)$,
$C(\psi)$ is a constant depending on $f$ and $\psi$, and $(u_n)$ is some
sequence
decaying to zero with rate determined by $f$ and $\mathcal{R}_\varepsilon
(\psi)$.  Notice that we make no assumption on the regularity of the
observable $\varphi$; when discussing the regularity class of observables, we
shall always be referring to the choice of the function $\psi$.  (This is not
atypical, although some existing results do require that both functions have
some minimum regularity.)

For brevity, we shall simply give an estimate for $u_n$ in the
statement of each result.
For each example we also have a Central Limit Theorem for those observables
which give summable decay of correlations, and are not coboundaries.  We
recall that a real-valued observable $\psi$ \emph{satisfies the Central Limit 
Theorem} for $f$ if there exists $\sigma>0$ such that for every interval
$J \subset \mathbb{R}$,
\[ \nu \left\{ x \in X : \frac{1}{\sqrt n} \sum_{j=0}^{n-1} \left(
   \varphi (f^j (x)) - \int \varphi d\nu \right) \in J \right\}
   \rightarrow \frac{1}{\sigma \sqrt{2 \pi}} \int_J e^{-\frac{t^2}
   {2 \sigma^2}} dt .\]

Note that the range of examples given in the following subsections is meant to be
illustrative rather than exhaustive, and so we shall miss out some simple
generalisations for which essentially the same results hold.  We shall instead
try to make clear the conditions needed to apply these results, and direct the
reader to the papers mentioned below for further examples which satisfy these
conditions.

\label{examples}

\subsection{Uniformly expanding maps}
\label{uniform}
Let $f:M \circlearrowleft$ be a $C^2$-diffeomorphism of a compact
Riemannian manifold.  We say $f$ is \emph{uniformly expanding} if
there exists $\lambda>1$ such that $\| Df_x v \| \geq \lambda \| v \|$ for
all $x \in M$, and all tangent vectors $v$.
Such a map admits an absolutely continuous invariant probability
measure $\mu$, which is unique and mixing.

\begin{theorem}
\label{uniformthm}

Let $\varphi \in \mathcal{L}^\infty (M,\mu)$, and let $\psi:M \rightarrow
\mathbb{R}$ be continuous.  Upper bounds are given
for $(u_n)$ as follows:
\begin{itemize}
\item if $\psi \in (R1)$, then $u_n = \mathcal{O} (\theta^n)$ for some
   $\theta \in (0,1)$;
\item if $\psi \in (R2,\gamma)$, for some $\gamma \in (0,1)$, then
$u_n = \mathcal{O} (e^{-n^{\gamma'}})$ for every $\gamma'<\gamma$;
\item if $\psi \in (R3,\gamma)$, for some $\gamma>1$, then
$u_n = \mathcal{O} (e^{-(\log n)^{\gamma'}})$ for every $\gamma'<\gamma$;
\item for any constant $C_\infty>0$ there exists $\zeta<1$
such that  if $\psi \in (R4,\gamma)$ for some $\gamma> \zeta^{-1}$, and
$\mathcal{R}_\infty (\psi) < C_\infty$, then
$u_n = \mathcal{O} (n^{1 - \zeta \gamma})$.

\end{itemize}

Furthermore, the Central Limit Theorem holds when $\psi \in (R4,\gamma)$
for sufficiently large $\gamma$, depending on $\mathcal{R}_\infty(\psi)$.

\end{theorem}

Such maps are generally regarded as being well understood, and in particular,
results of exponential decay of correlations for observables in $(R1)$
go back to the seventies, and the work of
Sinai, Ruelle and Bowen (\cite{Si}, \cite{R}, \cite{B}).
For a more modern perspective, see for instance the books of Baladi
(\cite{Ba}) and Viana (\cite{V2}).

I have not seen explicit claims of similar results for observables in
classes $(R2-4)$.  However, it is well known that any such map can be coded
by a one-sided full shift on finitely many symbols, so an analogous result
on shift spaces would be sufficient, and may well already exist.  The
estimates here are probably not sharp, particularly in the $(R4)$ case.

The other examples we consider are not in general reducible to finite
alphabet shift maps, so we can be more confident that the next set of
results are new.

\subsection{Maps with indifferent fixed points}
\label{indifferent}

These are perhaps the simplest examples of strictly non-uniformly expanding
systems.  Purely for simplicity, we restrict to the well known
case of the Manneville-Pomeau map.

\begin{theorem}
\label{indifferentthm}
\label{indiffthm}
Let $f:[0,1] \circlearrowleft$ be the map $f(x)=x + x^{1 + \alpha}$
($\mathrm{mod}$ 1), for some $\alpha \in (0,1)$, and let $\nu$ be the
unique acip for this system.
For $\varphi, \psi:[0,1] \rightarrow \mathbb{R}$ with $\varphi$ bounded and
$\psi$ continuous, for every constant $C_\infty>0$
there exists $\zeta<1$ such that if $\psi \in (R4,
\gamma)$ for some $\gamma> 2 \zeta^{-1}$, with $\mathcal{R}_\infty (\psi)<
C_\infty$, then
\begin{itemize}
\item if $\gamma=\zeta^{-1} (\tau+1)$, then $u_n = \mathcal{O}
   (n^{1 - \tau} \log n)$;
\item otherwise, $u_n = \mathcal{O} (\max (n^{1 - \tau},
   n^{2 - \zeta \gamma}))$;
\end{itemize}
where $\tau = \alpha^{-1}$.
In particular, when $\gamma>\frac{3}{\zeta}$ the Central Limit Theorem
holds.
\end{theorem}

In the case where $\psi \in (R4,\gamma)$ for \emph{every} large
$\gamma$, this gives $u_n = \mathcal{O}(n^{1-\frac{1}{\alpha}})$, which is
the bound obtained in \cite{Y2} for $\psi \in (R1)$.  We do not give separate
estimates for observables in classes $(R2)$ and $(R3)$, as we obtain
the same upper bound in each case.  Note that the polynomial
upper bound for $(R1)$ observables is known to be sharp (\cite{Hu}),
and hence the above gives a sharp bound in the $(R2)$ and $(R3)$ cases,
and for $(R4,\gamma)$ when $\gamma$ is large.

The above results apply in the more general 1-dimensional case considered in \cite
{Y2}, where in particular a finite number of expanding branches are allowed, and
it is assumed that $x f''(x) \approx x^\alpha$ near the indifferent fixed point.

In our remaining examples, estimates will invariably correspond to either
the above form, or that of Theorem \ref{uniformthm}, and we shall simply say which
is the case, specifying the parameter $\tau$ as appropriate.

\subsection{One-dimensional maps with critical points}
\label{BLvSintro}

Let us consider the systems of \cite{BLvS}.  These are
one-dimensional multimodal maps, where there is some
long-term growth of derivative along the critical orbits.
Let $f:I \rightarrow I$
be a $C^3$ interval or circle map with a finite
critical set $\mathcal{C}$ and no stable or neutral periodic orbit.  We
assume all critical points
have the same \emph{critical order} $l \in (1, \infty)$; this
means that for each $c \in \mathcal{C}$, there is some
neighbourhood in which $f$ can be written in the form
\[ f(x) = \pm | \varphi (x-c)|^l + f(c) \]
for some diffeomorphism $\varphi: \mathbb{R} \rightarrow \mathbb{R}$ fixing
$0$, with the $\pm$ allowed to depend on the sign of $x-c$.

For $c \in \mathcal{C}$, let $D_n(c) = |(f^n)' (f(c))|$.
From \cite{BLvS} we know there exists an acip $\mu$ provided
\[ \sum_n D_n^{-\frac{1}{2l-1}} (c) < \infty \, \, \, \forall c \in
\mathcal{C} . \]
If $f$ is not renormalisable on the support of $\mu$ then $\mu$ is mixing.

\begin{theorem}
\label{BLvSthm}

Let $\varphi \in \mathcal{L}^\infty (I,\mu)$, and let $\psi$ be continuous.

\emph{Case 1:} Suppose there exist $C>0$, $\lambda>1$ such that $D_n(c)
\geq C \lambda^n$ for all $n \geq 1$, $c \in \mathcal{C}$.  Then
we have estimates for $(u_n)$ exactly as in the uniformly expanding case
(Theorem \ref{uniformthm}).

\emph{Case 2:} Suppose there exist $C>0$, $\alpha>2l-1$ such that
$D_n (c) \geq C n^{\alpha}$ for all $n \geq 1$, $c \in \mathcal{C}$.
Then we have estimates for $(u_n)$ as in the indifferent fixed point case
(Theorem \ref{indiffthm}) for \emph{every} $\tau< \frac{\alpha-1}{l-1}$.

In particular, the Central Limit Theorem holds in either case when $\psi \in
(R4,\gamma)$ 
for sufficiently large $\gamma$, depending on
$\mathcal{R}_\infty(\psi)$.

\end{theorem}

Again, we have restricted our attention to some particular cases; analogous
results should be possible for the intermediate cases considered in \cite{BLvS}.
In particular, for the class of \emph{Fibonacci maps}
with quadratic critical points (see \cite{LM}) we obtain estimates as in
Theorem \ref{indiffthm} for \emph{every} $\tau>1$.

\subsection{Viana maps}
\label{vianamaps}

Next we consider the class of \emph{Viana maps}, introduced in \cite{V}.
These are examples of non-uniformly expanding maps in more than one
dimension, with sub-exponential decay of correlations for H\"older
observables.  They are notable for being possibly the first examples of
non-uniformly expanding systems in more than one dimension which admit an
acip, and also because the attractor, and many of its statistical properties,
persist in a $C^3$ neighbourhood of systems.

Let $a_0$ be some real number in $(1,2)$
for which $x=0$ is pre-periodic for the system $x \mapsto a_0 - x^2$.
We define a skew product $\hat{f}:S^1 \times
\mathbb{R} \circlearrowleft$ by
\[ \hat{f}(s,x)= (ds \textrm{ mod } 1, a_0 + \alpha \sin (2 \pi s)
- x^2), \]
where $d$ is an integer $\geq 16$, and $\alpha>0$ is a constant.
When $\alpha$ is sufficiently small, there is a compact interval $I
\subset (-2,2)$ for which $S^1 \times I$ is mapped strictly inside its own
interior, and $\hat{f}$ admits a unique acip, which is mixing for some
iterate, and has two positive Lyapunov exponents (\cite{V},\cite{AV}).
The same is also true for any $f$ in a sufficiently small $C^3$
neighbourhood $\mathcal{N}$ of $\hat{f}$.

Let us fix some small $\alpha$, and let $\mathcal{N}$ be a sufficiently
small
$C^3$ neighbourhood of $\hat{f}$ such that for every $f \in \mathcal{N}$ the
above properties hold.  Choose some $f \in \mathcal{N}$; if
$f$ is not mixing, we consider instead the
first mixing power.

\begin{theorem}
\label{vianathm}
For $\varphi \in \mathcal{L}^\infty (S^1 \times \mathbb{R},
\nu)$, $\psi \in (R4,\gamma)$,
we have estimates for $(u_n)$ as in the indifferent fixed point
case (Theorem \ref{indiffthm}) for \emph{every} $\tau>1$.

The Central Limit Theorem holds for $\psi \in (R4,\gamma)$
when $\gamma$ is sufficiently large, depending on $\mathcal{R}_\infty
(\psi)$.
\end{theorem}

Another way of saying the above is that if $\psi \in (R4,\gamma)$,
then $u_n = \mathcal{O} (n^{2 - \zeta \gamma})$, with the usual
dependency of $\zeta$ on $\mathcal{R}_\infty (\psi)$.
Note that for observables in $\cap_{\gamma>1} (R4,\gamma)$,
we get super-polynomial decay
of correlations, the same estimate as we obtain for H\"older observables
(though Baladi and Gou\"ezel have recently announced a stretched exponential
bound for H\"older observables - see {BG}).

There are a number of generalisations we could consider, such as
allowing $D \geq 2$ (\cite{BST} - note they require $f$ to be $C^\infty$
close to $\hat{f}$), or replacing $\sin (2 \pi s)$ by an arbitrary Morse
function.

\subsection{Non-uniformly expanding maps}

Finally, we discuss probably the most general context in which our
methods can currently be applied, the setting of \cite{ALP1}.  In
particular, this setting generalises that of Viana maps.

Let $f:M \rightarrow M$ be a transitive $C^2$ local diffeomorphism away from
a singular/critical set $\mathcal{S}$, with $M$ a compact
finite-dimensional Riemannian manifold.  Let $\mathrm{Leb}$ be a normalised
Riemannian volume form on $M$, which we shall refer to as Lebesgue measure,
and $d$ a Riemannian metric.  We assume $f$ is non-uniformly expanding, or
more precisely, there exists $\lambda>0$ such that
\begin{eqnarray} 
  \label{alp*1}
 \liminf_{n \rightarrow \infty} \frac{1}{n} \sum_{i=0}^{n-1}
   \log \| Df^{-1}_{f^i (x)} \|^{-1} \geq \lambda >0 .\end{eqnarray}

For almost every $x$ in $M$, we may define
\[ \mathcal{E} (x) = \min \left\{ N : \frac{1}{n} \sum_{i=0}^{n-1} \log \|
   Df^{-1}_{f^i (x)} \|^{-1} \geq \lambda/2 , \, \, \forall n
   \geq N \right\} .\]

The decay rate of the sequence $\mathrm{Leb} \{ \mathcal{E}(x)>n \}$ may
be considered to give a degree of hyperbolicity.  Where $\mathcal{S}$ is
non-empty, we need the following further assumptions, firstly on the
critical set.  We assume $\mathcal{C}$ is \emph{non-degenerate}, that is,
$m(\mathcal{C})=0$, and $\exists \beta>0$ such that $\forall x \in M \
\mathcal{C}$ we have $d(x,\mathcal{C})^\beta \lesssim \|Df_x v\| / \| v \|
\lesssim d(x,\mathcal{C})^{-\beta} \forall v \in T_x M$, and the functions
$\log \det Df$ and $\log \| Df^{-1} \|$ are locally Lipschitz with Lipschitz
constant $\lesssim d(x,\mathcal{C})^{-\beta}$.

Now let $d_\delta(x,\mathcal{S}) = d(x,\mathcal{S})$ when this is $\leq \delta$,
and $1$ otherwise.  We assume that for any $\varepsilon>0$ there exists
$\delta>0$ such that for Lebesgue a.e. $x \in M$,
\begin{eqnarray}
\label{alp*2}
 \limsup_{n \rightarrow \infty} \frac{1}{n} \sum_{j=0}^{n-1} - \log d_\delta
   (f^j(x),\mathcal{S}) \leq \varepsilon . \end{eqnarray}
We define a \emph{recurrence time}
\[ \mathcal{T}(x) = \min \left\{ N \geq 1 : \frac{1}{n} \sum_{i=0}^{n-1} -
  \log d_\delta (f^j (x),\mathcal{S}) \leq 2 \varepsilon , \, \,
  \forall n \geq N \right\} . \]
Let $f$ be a map satisfying the above conditions, and for which 
there exists $\alpha>1$ such that
\[ \mathrm{Leb} ( \{ \mathcal{E} (x)>n \textrm { or } \mathcal{T}(x)>n \}) =
   \mathcal{O} (n^{-\alpha}) .\]
Then $f$ admits an acip $\nu$ with respect to Lebesgue measure, and we may
assume $\nu$ to be mixing by taking a suitable power of $f$.

\begin{theorem}
\label{ALPthm}
For $\varphi \in \mathcal{L}^\infty (M,\nu), \psi \in (R4,\gamma)$, we have
estimates for $(u_n)$ as in the indifferent fixed point case (Theorem
\ref{indiffthm}), for $\tau=\alpha$.
Furthermore, when $\alpha>2$, the Central Limit Theorem holds for
$\psi \in (R4,\gamma)$ when $\gamma$ is sufficiently large for given
$\mathcal{R}_\infty (\psi)$.
\end{theorem}

\section{Young's tower}

\label{towers}

In the previous section, we indicated the variety of systems we may consider.
We shall now state the main technical result, and with it the conditions a
system must satisfy in order for our result to be applicable.  As verifying
that a system satisfies such conditions is often considerable work, we refer
the reader to those papers mentioned in each of the previous subsections for
full details.

The relevant setting for our arguments will be the \emph{tower} object
introduced by Young in \cite{Y2}, and we recap its definition.
We start with a map $F^R:(\Delta_0,m_0) \circlearrowleft$, where
$(\Delta_0,m_0)$ is a finite measure space.  This shall represent the
\emph{base} of the tower.  We assume there exists a partition (mod $0$)
$\mathcal{P} = \{ \Delta_{0,i}: i \in \mathbb{N} \}$ of $\Delta_0$, such
that $F^R|\Delta_{0,i}$ is an injection onto $\Delta_0$ for each
$\Delta_{0,i}$.  We require that the partition generates, i.e. that
$\bigvee_{j=0}^\infty (F^R)^{-j} \mathcal{P}$ is the trivial partition into
points.  We also choose a return time function $R:\Delta_0 \rightarrow
\mathbb{N}$, which must be constant on each $\Delta_{0,i}$.

We define a tower to be any map $F:(\Delta,m) \circlearrowleft$ determined by
some $F^R$, $\mathcal{P}$, and $R$ as follows.  Let $\Delta =
\{ (z,l) : z \in \Delta_0, \, l<R(z) \}$.  For convenience let $\Delta_l$
refer to the set of points $(\cdot , l)$ in $\Delta$.
This shall be thought of as the $l$th \emph{level} of $\Delta$.  (We shall
freely confuse the \emph{zeroth level} $\{ (z,0): z \in
\Delta_0 \} \subset \Delta$ with
$\Delta_0$ itself.  We shall also happily refer to points in $\Delta$ by a
single letter $x$, say.)  We write $\Delta_{l,i}= \{ (z,l): z \in
\Delta_{0,i} \}$ for $l<R(\Delta_{0,i})$.  The partition of $\Delta$ into the
sets $\Delta_{l,i}$ shall be denoted by $\eta$.

The map $F$ is then defined as follows:
\[ F(z,l) = \left\{ \begin{array}{rcl} (z,l+1)&& \textrm{if } l+1<R(z) \\
   (F^R(z),0)&& \textrm{otherwise.} \end{array} \right. \]

We notice that the map $F^{R(x)}(x)$ on $\Delta_0$ is identical to $F^R(x)$,
justifying our choice of notation.  Finally, we define a notion of
\emph{separation time}; for $x,y \in \Delta_0$, $s(x,y)$ is defined to be
the least integer $n \geq 0$ s.t. $(F^R)^nx,(F^R)^ny$ are in different
elements of $\mathcal{P}$.  For $x,y \in$ some $\Delta_{l,i}$, where
$x=(x_0,l)$, $y=(y_0,l)$, we set
$s(x,y):=s(x_0,y_0)$; for $x,y$ in different elements of
$\eta$, $s(x,y)=0$.

We say that the \emph{Jacobian} $JF^R$ of $F^R$ with respect
to $m_0$ is the real-valued function such that for any measurable set $E$
on which $JF^R$ is injective,
\[ m_0 (F^R (E)) = \int_E JF^R dm_0 .\]
We assume $JF^R$ is uniquely defined, positive, and finite $m_0$-a.e.
We require some further assumptions.

\begin{itemize}
\item \emph{Measure structure:} Let $\mathcal{B}$ be the $\sigma$-algebra of
$m_0$-measurable sets.  We assume that all elements of $\mathcal{P}$ and each
$\bigvee_{i=0}^{n-1} (F^R)^{-i} \mathcal{P}$ belong to $\mathcal{B}$, and
that $F^R$ and $(F^R | \Delta_{0,i})^{-1}$ are measurable functions.  We then
extend $m_0$ to a measure $m$ on $\Delta$ as follows: for $E \subset
\Delta_l$, any $l \geq 0$, we let
$m(E)=m_0(F^{-l}E)$, provided that $F^{-l}E \in \mathcal{B}$.
Throughout, we shall assume that any sets we choose are measurable.  Also,
whenever we say we are choosing an arbitrary point $x$, we shall assume it is
a \emph{good} point, i.e. that each element of its orbit is contained within
a single element of the partition $\eta$, and that $JF^R$ is well-defined
and positive at each of these points.

\item \emph{Bounded distortion:} There exist $C>0$ and $\beta<1$ s.t. for
$x,y \in$ any $\Delta_{0,i} \in \mathcal{P}$,
\[ \left| \frac{JF^R(x)}{JF^R(y)} -1 \right| \leq C \beta^{s(F^R x,F^R y)} .\]

\item \emph{Aperiodicity:} We assume that $\mathrm{gcd} \{ R(x): x \in
\Delta_0 \}=1$.  This is a necessary and sufficient condition for mixing (in
fact, for exactness).

\item \emph{Finiteness:} We assume $\int R dm_0 < \infty$.  This tells us
that $m(\Delta)<\infty$.

\end{itemize}

Let $F:(\Delta,m) \circlearrowleft$ be a tower, as defined above.  We
define classes of observable similar to those we consider on the manifold,
but characterised instead in terms of the separation time $s$ on $\Delta$.
Given a bounded function $\psi:\Delta \rightarrow \mathbb{R}$, we define
the \emph{variation} for $n \geq 0$:

\[ v_n (\psi) = \sup \{ | \psi(x) - \psi(y)|: s(x,y) \geq n \}. \]
Let us use this to define some regularity classes:

\emph{Exponential case:} $\psi \in (V1, \gamma)$, $\gamma \in (0,1)$, if
$v_n(\psi) = \mathcal{O} (\gamma^n)$;

\emph{Stretched exponential case:} $\psi \in (V2,\gamma)$, $\gamma \in
(0,1)$, if $v_n(\psi) = \mathcal{O}(\exp \{ -n^\gamma \})$;

\emph{Intermediate case:} $\psi \in (V3,\gamma)$, $\gamma>1$, if
$v_n(\psi) = \mathcal{O} (\exp \{ -(\log n)^\gamma \} )$;

\emph{Polynomial case:} $\psi \in (V4,\gamma)$, $\gamma>1$, if $v_n(\psi)
= \mathcal{O} (n^{-\gamma})$.

We shall see that the classes (V1-4) of regularity correspond naturally
with the classes (R1-4) of regularity on the manifold respectively, under
fairly weak assumptions on the relation between the system and the tower
we construct for it.  (We shall discuss this further in \S
\ref{applications}.)  These classes are essentially those defined in
\cite{P}, although there the functions are considered to be \emph{potentials}
rather than observables.

We now state the main technical result.

\begin{theorem}
\label{Ythm2}
\label{towerthm}
Let $F:(\Delta,m) \circlearrowleft$ be a tower satisfying the assumptions
stated above.  Then $F:(\Delta,m)\circlearrowleft$ admits a
unique acip $\nu$, which is mixing.  Furthermore, for all $\varphi,\psi
\in \mathcal{L}^\infty (\Delta,m)$,
\[ \left| \int (\varphi \circ F^n ) \psi d\nu - \int \varphi d\nu \int \psi
   d\nu \right| \leq \| \varphi \|_\infty C(\psi) u_n ,\]
where $C(\psi)>0$ is some constant, and $(u_n)$ is a sequence converging
to zero at some rate determined by $F$ and $v_n(\psi)$.  In particular:

\emph{Case 1:} Suppose $m_0 \{ R>n \} = \mathcal{O} (\theta^n)$, some $\theta
\in (0,1)$.  Then
\begin{itemize}
\item if $\psi \in (V1,\gamma)$ for some $\gamma
\in (0,1)$, then
$u_n = \mathcal{O}(\theta^n)$ for some $\theta \in (0,1)$;
\item if $\psi \in (V2,\gamma)$ for some $\gamma
\in (0,1)$, then
$u_n = \mathcal{O}(e^{-n^{\gamma'}})$ for every $\gamma'<\gamma$;
\item if $\psi \in (V3,\gamma)$ for some
$\gamma>1$, then
$u_n = \mathcal{O}(e^{-(\log n)^{\gamma'}})$ for every $\gamma'<\gamma$;
\item for any constant $C_\infty>0$, there exists $\zeta<1$ such that
if $\psi \in (V4,\gamma)$ for some $\gamma>\frac{1}{\zeta}$, and $v_0
(\psi)< C_\infty$, then $u_n= \mathcal{O}(n^{1-\zeta \gamma})$.
\end{itemize}

\emph{Case 2:} Suppose $m_0 \{ R>n \} = \mathcal{O} (n^{-\alpha})$ for some
$\alpha>1$.  Then for every $C_\infty>0$ there exists $\zeta<1$ 
such that if $\psi \in (V4,\gamma)$ for some
$\gamma>\frac{2}{\zeta}$, with $v_0(\psi) < C_\infty$, then
\begin{itemize}
\item if $\gamma = \frac{\alpha+1}{\zeta}$,
   $u_n=\mathcal{O}(n^{1-\alpha} \log n)$;
\item otherwise, $u_n=\mathcal{O}\left( \max (n^{1-\alpha},n^{2-\zeta
   \gamma}) \right)$.
\end{itemize}

\end{theorem}

The existence of a mixing acip is proved in \cite{Y2}, as is the result in
the case $\psi \in (V1)$.
As a corollary of the above, we get a Central Limit Theorem in the cases
where the rate of mixing is summable.

\begin{corr}
\label{clt}
Suppose $F$ satisfies the above assumptions, and $m_0 \{ R>n \} = 
\mathcal{O}(n^{-\alpha})$, for some $\alpha>2$.  Then the Central Limit
Theorem is satisfied for $\psi \in (R4,\gamma)$ when $\gamma$ is
sufficiently large, depending on $F$ and $v_0(\psi)$.
\end{corr}

In \S \ref{applications} we shall give the exact conditions needed on a
system in order to apply the above results.

\section{Overview of method}
\label{method}

Our strategy in proving the above theorem is to generalise a
\emph{coupling} method introduced by Young in \cite{Y2}.  Our argument
follows closely the line of approach of that paper, and we give an outline
of the key ideas here.

First, we need to reduce the problem to one in a slightly different
context.
Given a system
$F:(\Delta,m) \circlearrowleft$, we define a \emph{transfer operator}
$F_*$ which, for any measure $\lambda$ on $\Delta$ for which $F$ is
measurable, gives a measure $F_* \lambda$ on $\Delta$ defined by
\[ (F_* \lambda) (A) = \lambda (F^{-1} A) \]
whenever $A$ is a $\lambda$-measurable set.  Clearly any $F$-invariant
measure is a fixed point for this operator.  Also, a key property of $F_*$
is that for any function $\phi: \Delta \rightarrow \mathbb{R}$,
\[ \int \phi \circ F d\lambda = \int \phi d(F_* \lambda) .\]

Next, we define a \emph{variation norm} on $m$-absolutely continuous
\emph{signed measures}, that is, on the difference between any two (positive)
measures which are absolutely continuous.  Given two such measures
$\lambda,\lambda'$, we write
\[ | \lambda - \lambda'| := \int \left| \frac{d\lambda}{dm} - \frac{
   d\lambda'}{dm} \right| dm .\]
Now let us fix an acip $\nu$ and choose observables $\varphi \in
\mathcal{L}^\infty (\Delta,\nu)$, $\psi \in \mathcal{L}^1 (\Delta,\nu)$,
with $\inf \psi >0$, $\int \psi d\nu =1$.  We have
\begin{eqnarray*} \int (\varphi \circ F^n) \psi d\nu &=&
   \int (\varphi \circ F^n) d(\psi \nu) \\
   &=& \int \varphi d(F^n_* (\psi \nu)) ,\end{eqnarray*}
where $\psi \nu$ denotes the unique measure which has density $\psi$
with respect to $\nu$.  So
\begin{eqnarray*}
   \left| \int (\varphi \circ F^n) \psi d\nu - \int \varphi d\nu \int
       \psi d\nu \right| &\leq& \| \varphi \|_\infty \int
       \left| \frac{dF^n_* (\psi \nu)}{dm} - \frac{d\nu}{dm} \right| dm \\
    &=& \| \varphi \|_\infty |F^n_* (\psi \nu) - \nu| .\end{eqnarray*}

Hence we may reduce the problem to one of estimating the rate at which certain
measures converge to the invariant measure, in terms of the variation norm.
In fact, it will be useful to consider the more general question of
estimating $|F^n_* \lambda - F^n_* \lambda'|$ for a pair of measures
$\lambda$, $\lambda'$ whose densities with respect to $m$ are of some
given regularity.  (We shall require an estimate in the case $\lambda' \neq
\nu$ when we consider the Central Limit Theorem.)

Let us now outline the main argument.  We work with two copies of the system,
and the direct product $F \times F: (\Delta \times \Delta, m \times m)
\circlearrowleft$.  Let $P_0 = \lambda \times \lambda'$, and consider it to
be a measure on $\Delta \times \Delta$.  If we let $\pi, \pi': \Delta
\times \Delta \rightarrow \Delta$ be the projections onto the first and
second coordinates respectively, we have that
\begin{eqnarray*}
   |F^n_* \lambda - F^n_* \lambda' | &=& | \pi_* (F \times F)^n_* P_0
		- \pi_*' (F \times F)^n_* P_0| \\
	&\leq& 2 |(F \times F)^n_* P_0 | .\end{eqnarray*}

Our strategy will involve summing the differences between the two projections
over small regions of the space, only comparing them at convenient times
that vary with the region of space we are considering.  At each of these
times, we shall subtract some measure from both coordinates so that the
difference is unaffected, yet the total measure of $P_0$ is reduced, giving
an improved upper bound on the difference.

The key difference between our method and that of \cite{Y2} is that we
introduce a sequence $(\varepsilon_n)$, which shall represent the rate at
which we attempt to subtract measure from $P_0$.  When the densities of
$\lambda$, $\lambda'$ are of class $(V1)$, $(\varepsilon_n)$ can be taken
to be a small constant, and the method here reduces to that of \cite{Y2};
however, by allowing sequences $\varepsilon_n \rightarrow 0$, we may also
consider measure densities of weaker regularity.

We shall see that it is possible to define an induced map $\hat{F}:
\Delta \times \Delta \rightarrow \Delta_0 \times \Delta_0$ for which there
is a partition $\hat{\xi}_1$ of $\Delta \times \Delta$, with every element
mapping injectively onto $\Delta_0 \times \Delta_0$ under $\hat{F}$.  In
fact, there is a stopping time $T_1:\hat{\xi}_1 \rightarrow \mathbb{N}$ such
that for each $\Gamma \in \hat{\xi}_1$,
$\hat{F}|\Gamma = (F \times F)^{T_1 (\Gamma)}$.  If we choose some
$\Gamma \in \hat{\xi}_1$, then $\hat{F}_* (P_0 | \Gamma)$ is a measure on
$\Delta_0 \times \Delta_0$.  
The density of $P_0$ with respect to $m \times m$ has essentially the same
regularity as the measures $\lambda$, $\lambda'$, and the density of
$\hat{F}_* (P_0 | \Gamma)$ will be similar, except possibly weakened
slightly by any irregularity in the map $\hat{F}$.  (We shall see that the
map $\hat{F}$ is not too irregular.)

Let
\[ c(\Gamma) = \inf_{w \in \Delta_0 \times \Delta_0} \frac{d \hat{F}_*
   (P_0 | \Gamma)}{d(m \times m)} (w) .\]
For any $\varepsilon_1 \in [0,1]$, we may write
\[ \hat{F}_* (P_0 | \Gamma) = \varepsilon_1 c(\Gamma)
   (m \times m | \Delta_0 \times
   \Delta_0) + \hat{F}_* (P_1 | \Gamma) \]
for some (positive) measure $P_1|\Gamma$; this is uniquely defined since
$\hat{F}| \Gamma$ is injective.  Essentially, we are subtracting some
amount of mass from the measure $\hat{F}_* (P_0 | \Gamma)$.  Moreover, we
are subtracting it equally from both coordinates; this means that writing
$\Gamma = A \times B$ and $k = T_1 (\Gamma)$, the distance between the
measures $F^k_* (\lambda |A)$ and $F^k_* (\lambda'|B)$, both defined on
$\Delta_0$, is unaffected.  However, we also see that the remaining
measure $\hat{F}_* (P_1 | \Gamma)$ has smaller total mass, and this is an
upper bound for $|F^k_* (\lambda|A) - F^k_* (\lambda'|B)|$.

We fix an $\varepsilon_1$, and perform this subtraction of measure for
each $\Gamma \in \hat{\xi}_1$, obtaining a measure $P_1$ defined on $\Delta
\times \Delta$.  The total mass of $P_1$ represents the difference between
$F^n_* \lambda$ and $F^n_* \lambda'$ at time $n=T_1$, taking into account
that $T_1$ is not constant over $\Delta \times \Delta$.  Clearly, we
obtain the best upper bound by taking $\varepsilon_1=1$; however, we shall
see that it is to our advantage to choose some smaller value for
$\varepsilon_1$.

We choose a sequence $(\varepsilon_n)$, and proceed inductively as follows.
First, we define a sequence of partitions $\{ \hat{\xi}_i \}$ such that
$\Gamma \in \hat{\xi}_i$ is mapped injectively onto $\Delta_0 \times
\Delta_0$ under $\hat{F}^i$.  Now given the measure $P_{i-1}$, we take an
element $\Gamma \in \hat{\xi}_i$ and consider the measure
$\hat{F}^i_* (P_{i-1}|\Gamma)$ on $\Delta_0 \times \Delta_0$.  Here, we let
\[ c(\Gamma) = \inf_{w \in \Delta_0 \times \Delta_0}
   \frac{d\hat{F}^i_* (P_{i-1} | \Gamma)}
   {d(m \times m | \Delta_0 \times \Delta_0)} (w) ,\]
and specify $P_i|\Gamma$ by
\[ \hat{F}^i_* (P_{i-1} | \Gamma) = \varepsilon_i c(\Gamma) ( m \times m | \Delta_0
   \times \Delta_0) + \hat{F}^i_* (P_i | \Gamma) .\]

As before, we construct a measure $P_i$, the total mass of which gives an
upper bound at time $T_i$.  To fully determine the sequence $\{ P_i \}$,
it remains to choose a sequence $(\varepsilon_i)$.  Our choice relates to
the regularity of the densities $\frac{d\lambda}{dm},\frac{d\lambda'}{dm}$.
This is relevant because the method requires that the family of measure
densities
\[ \left\{ \frac{d \hat{F}^i_* (P_{i-1} | \Gamma)}{d(m \times m)} :
   i \geq 1, \Gamma \in \hat{\xi}_i \right\} \]
has some uniform regularity.  (In fact, we require that the $\log$ of each
of the above densities is suitably regular.)  We require this in order that
at the $i$th stage of the procedure, 
when we subtract an $\varepsilon_i$ proportion of the minimum local density,
this corresponds to a similarly large proportion of the average density.
Hence, provided this regularity is maintained, the total mass of $P_{i-1}$ is
decreased by a similar proportion.

When we subtract a constant from a density as above, this weakens the
regularity.
However, at the next step of the procedure, we work with elements of the
partition $\hat{\xi}_{i+1}$.  Since these sets are smaller, we regain some
regularity by working with measures on $\Delta_0 \times \Delta_0$ pushed
forward from such sets.  That is, we expect the densities
\[ \left\{ \frac{d\hat{F}^{i+1}_* (P_i | \Gamma)}{d(m \times m)} : \Gamma
   \in \hat{\xi}_{i+1} \right\} \]
to be more regular than the densities
\[ \left\{ \frac{d\hat{F}^i_* (P_i | \Gamma)}{d(m \times m)} : \Gamma \in
   \hat{\xi}_i \right\} .\]
(This relies on the map $\hat{F}$ being smooth enough that another
application of the operator $\hat{F}_*$ doesn't much affect the regularity.)

The degree of regularity we gain in this way depends on the initial
regularity of $\Phi$, and hence of
$\frac{d\lambda}{dm},\frac{d\lambda'}{dm}$, with respect to the sequence
of partitions.  In the usual case, where $\frac{d\lambda}{dm},
\frac{d\lambda'}{dm} \in (V1)$, the regularity we gain
each time we refine the partition is similar to the regularity
we lose when we subtract a small
constant proportion of the density; hence we may take every $\varepsilon_i$
to be a small constant $\varepsilon$.  Where the initial regularities are
not so good, we gain less regularity from refining the partition, and so
we may only subtract correspondingly less measure.

For this reason, outside the $(V1)$ case we shall require that the sequence
$(\varepsilon_i)$ converges to zero at some minimum rate.  However, if
$(\varepsilon_i)$ decays faster than necessary, we will simply obtain a
suboptimal bound.  So part of the problem is to try to choose a sequence
$(\varepsilon_i)$ decaying as slowly as is permissible.
We shall also need to take into account the stopping time $T_1$ (which is
unbounded), in order to estimate the speed of convergence in terms of the
original map $F$.

\section{Related work}

\label{related}
Let us now mention some other results concerning estimates on
decay of correlations for non-H\"older observables.  Most of these
are stated in the context of one-sided finite alphabet shift maps,
or subshifts of finite type.  (For a comprehensive discussion of
shift maps and equilibrium measures, we suggest the book of Baladi,
\cite{Ba}.)
Shift maps are relatively simple dynamical systems, but are often
used to \emph{code} more complicated systems via a
semi-conjugacy, in much the same way that each
of the examples we consider can be represented by a suitable
\emph{tower} (see \ref{applications}).
Where a system $F:X \rightarrow X$
being coded has an invariant measure $\mu$ which is absolutely
continuous with respect to Lebesgue measure, $\mu$ is an
equilibrium measure for the potential $\phi = - \log JF$,
where $JF$ is the Jacobian with respect to Lebesgue measure.
Most results for shift maps work with an equilibruim measure
given by a potential $\phi$ which is H\"older continuous (in
terms of the usual metric on shift spaces - two sequences are
said to be distance $2^{-n}$ apart if they agree for exactly
the first $n$ symbols).
This assumption corresponds to assuming good distortion for $JF$.  

The have been various results concerned primarily with weakening
the assumption on the regularity of $\phi$, and obtaining (slower)
upper bounds for the rate of mixing with respect to the corresponding
equilibrium measures.
Kondah, Maume and Schmitt (\cite{KMS}) used a method of
Birkhoff cones and projective metrics, Bressaud,
Fernandez and Galves (\cite{BFG}) used a coupling method
(different from the one described here), with estimates given
in terms of \emph{chains of complete connections}, and 
Pollicott (\cite{P}) introduced a method involving
composing transfer operators with conditional expectations.
Each of these results has slightly different assumptions
and gives slightly different estimates, but in each case
a number of different classes of potentials are considered, 
and estimates are given for for observables of some similar
regularity to (usually \emph{not much worse} than) the
potential.  In particular, in all three examples polynomial
mixing is given for a potential and observables with
variations decaying at suitable polynomial rates.

In addition, Fisher and Lopes (\cite {FL}) and Isola (\cite{I})
have obtained polynomial decay of correlations for some
specific classes of potentials on the full 2-shift, each for a
class of observables not unlike our $(V4)$ class.

We emphasise that each of the above results concerns only shifts
(or subshifts) on a \emph{finite alphabet}.  Furthermore, with the
exception of uniformly expanding maps (\S \ref{uniform}), none of
the examples we consider may be coded (with a sufficiently
regular semi-conjugacy)  by shifts with finite alphabets, due to
distortion considerations.  Hence outside the uniform case, there
is no direct overlap between the above results and
our results. These are essentially generalisations in
a different direction,  but are worth mentioning if only
to note the variety of different methods that have been
applied to obtain results for non-H\"older observables in
slightly different contexts.

We now move on from shift spaces to mention a result
on towers.  Following the method of \cite{KMS} mentioned
above, Maume-Deschamps (\cite{M}) essentially reproduced
Young's results (\cite{Y2}), with some slightly weaker
estimates.  Together with Buzzi (\cite{BM}), the method has
been extended to allow the bounded distortion assumption of
\S \ref{towers} to be replaced by a condition on the variation
of $JF$, and also to allow observables of similar regularity.
One important difference is that the variation is defined in
terms of $JF$ and a partition of the whole tower (the
partition $\eta$ of \S \ref{towers}) rather than $JF^R$ and
the base partition, as here.  This significantly reduces
the classes of observables they can consider; for instance,
for any tower with unbounded return times, our class $(V1)$
contains many functions that cannot be dealt with at all
by their method.  It is not clear that any estimates can be
obtained for our examples, except for the uniformly
expanding case; for Manneville-Pomeau maps, for instance
(see \S \ref{indifferent}), while we can construct the
same tower in their context, the semi-conjugacy between the
tower and the
system is not regular enough to give any comparable
results.  Some applications are given in the related paper
\cite{BM2}, including certain multi-dimensional piecewise
expanding affine maps.

We also note that a similar result to that of \cite{BM} was
obtained by Holland (\cite{H1}), using a coupling method very
similar to the one we use here.

Finally, we mention a result which applies directly to certain
non-uniformly expanding systems, rather than to a symbolic
space or tower.  Pollicott and Yuri (\cite{PY}) consider a class
of maps of
arbitrary dimension with a single indifferent periodic orbit
and a given Markov structure, including in particular the
Manneville-Pomeau interval maps.  The class of observables
considered is dynamically defined; each observable is required
to be Lipschitz with respect to a Markov partition corresponding
to some induced map, chosen to have good distortion properties.
This class includes all functions which are Lipschitz with
respect to the manifold, and while some estimates are weaker
than comparable results for H\"older observables, bounds are
obtained for some observables which cannot be dealt with at all
by our methods, such as certain unbounded functions.

\section{Coupling}
\label{coupling}

Over the next few sections, we give the proof of the main technical theorem.

Let $F:(\Delta,m) \circlearrowleft$ be a tower, as defined in \S \ref{towers}.
Let $\mathcal{I}= \{ \varphi : \Delta \rightarrow \mathbb{R} \, | \, v_n
(\varphi) \rightarrow 0 \}$, and
let $\mathcal{I}^+ = \{ \varphi \in \mathcal{I} : \inf \varphi >0 \}$.
We shall work with probability measures whose densities with respect to
$m$ belong to $\mathcal{I}^+$.  We see
\begin{eqnarray}
\label{cpeq}
 \left| \frac{\varphi(x)}{\varphi(y)} -1 \right| = \frac{1}{\varphi(y)}
   | \varphi(x) - \varphi(y)| \leq C_\varphi v_{s(x,y)} (\varphi) ,
\end{eqnarray}
where $C_\varphi$ depends on $\inf \varphi$.

Let $\lambda$, $\lambda'$ be measures with $\frac{d\lambda}{dm}$,
$\frac{d\lambda'}{dm} \in \mathcal{I}^+$, and let
$P=\lambda \times \lambda'$.
For convenience, we shall write $v_n(\lambda)
=v_n(\frac{d\lambda}{dm})$ for such measures $\lambda$, and
use the two notations interchangeably.  We let $C_\lambda = C_\varphi$ above,
where $\varphi = \frac{d\lambda}{dm}$.  We shall write $\nu$ for the unique
acip for $F$, which is equivalent to $m$.

We consider the direct product $F \times F : (\Delta \times \Delta, m \times
m) \circlearrowleft$, and specify a return function to $\Delta_0 \times
\Delta_0$.  We first fix $n_0>0$ to
be some integer large enough that $m(F^{-n}\Delta_0 \cap \Delta_0) \geq$ some
$c>0$ for all $n \geq n_0$.  Such an integer exists since
$\nu$ is mixing and equivalent to $m$.  Now we let $\hat{R}(x)$ be the
first arrival time to $\Delta_0$ (setting $\hat{R}|\Delta_0 \equiv 0)$.
We define a sequence
$\{ \tau_i \}$ of stopping time functions on $\Delta \times \Delta$ as
follows:
\begin{eqnarray*}
   \tau_1(x,y)&=&n_0 + \hat{R}(F^{n_0}x), \\
   \tau_2(x,y)&=&\tau_1 + n_0 + \hat{R}(F^{\tau_1 + n_0} y), \\
   \tau_3(x,y)&=&\tau_2 + n_0 + \hat{R}(F^{\tau_2 + n_0} x),
\end{eqnarray*}
and so on, alternating between the two coordinates $x,y$ each time.
Correspondingly, we shall define an increasing sequence
$\xi_1 < \xi_2 < \ldots$ of partitions of
$\Delta \times \Delta$, according to each $\tau_i$.  First, let
$\pi$, $\pi'$ be the coordinate projections of $\Delta \times \Delta$ onto
$\Delta$, that is, $\pi(x,y):=x$, $\pi'(x,y):=y$.  At each stage we refine the
partition according to one of the two coordinates, alternating between the two
copies of $\Delta$.  First, $\xi_1$ is given by taking the partition into
rectangles $E \times \Delta$, $E \in \eta$, and refining so that $\tau_1$ is
constant on each element $\Gamma \in \xi_1$, and $F^{\tau_1}|\pi(\Gamma)$ is
for each $\Gamma$ an injection onto $\Delta_0$.  To be precise, we write
\[ \xi_1(x,y)= \left( \bigvee_{j=0}^{\tau_1(x,y)-1} F^{-j}\eta \right)(x)
   \times \Delta,\]
using throughout the convention that for a partition $\xi$, $\xi(x)$ denotes
the element of $\xi$ containing $x$.
Subsequently, we say $\xi_i$ is the refinement of $\xi_{i-1}$ such that each
element of $\xi_{i-1}$ is partitioned in the first (resp. second) coordinate
for $i$ odd (resp. even) so that $\tau_i$ is constant on each element $\Gamma
\in \xi_i$, and $F^{\tau_i}$ maps $\pi(\Gamma)$ (resp. $\pi'(\Gamma)$)
injectively onto $\Delta_0$.

We define $T$ to be the smallest $\tau_i$, $i \geq 2$, with
$(F \times F)^{\tau_i} (x,y) \in \Delta_0 \times \Delta_0$.  This is
well-defined $m$-a.e. since $\nu \times \nu$ is ergodic (in
fact, mixing).  Note that this is not necessarily the first
return time to $\Delta_0 \times \Delta_0$ for $F \times F$.  We now consider
the simultaneous
return function $\hat{F}:=(F \times F)^T$, and partition $\Delta \times
\Delta$ into regions which $\hat{F}^n$ maps injectively onto $\Delta_0 \times
\Delta_0$.

For $i \geq 1$ we let $T_i$ be the time
corresponding to the $i$th iterate of $\hat{F}$, i.e. $T_1 \equiv T$, and
for $i \geq 2$, 
\[ T_i (z) = T_{i-1} (z) + T(\hat{F}^{i-1} z) .\]
Corresponding to $\{ T_i \}$ we define a sequence of partitions 
$\eta \times \eta \leq \hat{\xi}_1 \leq \hat{\xi}_2 \leq \ldots$ of $\Delta
\times \Delta$ similarly to before, such that for each $\Gamma \in
\hat{\xi}_n$, $T_n | \Gamma$ is constant and $\hat{F}^n$ maps $\Gamma$
injectively onto $\Delta_0 \times \Delta_0$.  It will be convenient to define
a separation time $\hat{s}$ with respect to $\hat{\xi}_1$; $\hat{s}(w,z)$ is
the smallest
$n \geq 0$ s.t. $\hat{F}^n w$, $\hat{F}^n z$ are in different elements of
$\hat{\xi}_1$.  We notice that if $w=(x,x')$, $z=(y,y')$, then
$\hat{s}(w,z) \leq \min (s(x,y),s(x',y'))$.

Let $\varphi = \frac{d\lambda}{dm}$, $\varphi' =
\frac{d\lambda'}{dm}$, and let $\Phi = \frac{dP}{dm \times m} = \varphi \cdot
\varphi'$.
We first consider the regularity of $\hat{F}$ and $\Phi$ with respect to the
separation time $\hat{s}$.

\begin{sublemma}
\label{sbprodreg}
\begin{verbatim}

\end{verbatim}
\begin{enumerate}
\item For all $w,z \in \Delta \times \Delta$ with $\hat{s}(w,z) \geq n$,
\[ \left| \log \frac{J\hat{F}^n(w)}{J\hat{F}^n(z)} \right| \leq C_{\hat{F}}
   \beta^{\hat{s}(\hat{F}^n w,\hat{F}^n z)} \]
for $C_{\hat{F}}$ depending only on $F$.
\item For all $w,z \in \Delta \times \Delta$,
\[ \left| \log \frac{\Phi(w)}{\Phi(z)} \right| \leq C_\Phi v_{\hat{s}(w,z)}
   (\Phi) , \]
where $v_n(\Phi):= \max \left( v_n(\varphi),
v_n(\varphi') \right)$, and $C_\Phi = C_\varphi + C_{\varphi'}$,
where $C_\varphi$,$C_{\varphi'}$ are the constants given in (\ref{cpeq})
above, corresponding
to $\varphi,\varphi'$ respectively.
\end{enumerate}
\end{sublemma}

\emph{Proof:}
Let $w=(x,x')$, $z=(y,y')$.  When $\hat{s}(w,z) \geq n$, there exists $k \in
\mathbb{N}$ with $\hat{F}^n \equiv (F \times F)^k$ when restricted to the
element of $\hat{\xi}_n$ containing $w,z$.  So
\begin{eqnarray*}
 \left| \log \frac{J\hat{F}^n (w)}{J\hat{F}^n (z)} \right| &=&
   \left| \log \frac{JF^k(x) JF^k(x')}{JF^k(y) JF^k(y')} \right| \\
 &\leq&
   \left| \log \frac{JF^k(x)}{JF^k(y)} \right| + \left| \log
   \frac{JF^k(x')}{JF^k(y')} \right| .\end{eqnarray*}
Let $j$ be the number of times $F^i(x)$ enters $\Delta_0$, for $i=1,\ldots
,k$.  We have
\begin{eqnarray*}
  \left| \log \frac{JF^k (x)}{JF^k (y)} \right| &\leq& \sum_{i=1}^j
   C \beta^{s(F^k x,F^k y)+(j-i)} \\
  &\leq& C' \beta^{s(F^k x,F^k y)}
\end{eqnarray*}
for some $C'>0$, and similarly for $x',y'$.  So
\[ \left| \log \frac{J\hat{F}^n(w)}{J\hat{F}^n(z)}\right| \leq C_{\hat{F}}
   \beta^{\hat{s}(\hat{F}^n w, \hat{F}^n z)} \]
for some $C_{\hat{F}}>0$.  For the second part, 
we have
\begin{eqnarray*}
 \left| \log \frac{\Phi(w)}{\Phi(z)} \right| &\leq& \left| \log
   \frac{\varphi(x)}{\varphi(y)} \right| + \left| \log \frac{\varphi'(x')}
   {\varphi'(y')} \right| \\
 &\leq& C_\varphi v_{s(x,y)} (\varphi) + C_{\varphi'}
   v_{s(x',y')}(\varphi') \\
 &\leq& C_\Phi v_{\hat{s}(w,z)}(\Phi) .\end{eqnarray*}
\hfill $\square$

We now come to the core of the argument.  We choose a sequence
$(\varepsilon_i) < 1$, which represents the proportion of $P$ we
try to subtract at each step of the construction.  Let $\psi_0 \equiv
\tilde{\psi}_0 \equiv \Phi$.  We proceed as
follows; we \emph{push-forward} $\Phi$ by $\hat{F}$ to obtain the function
$\psi_1(z):= \frac{\Phi(z)}{J\hat{F}(z)}$.  On each element $\Gamma \in
\hat{\xi}_1$ we subtract the constant $\varepsilon_1\inf\{\psi_1(z) : z \in
\Gamma \}$ from the density $\psi_1 | \Gamma$.  We continue inductively,
pushing forward by dividing the density by $J\hat{F}(\hat{F}z)$ to get
$\psi_2(z)$, subtracting $\varepsilon_2 \inf\{\psi_2(z) :  z \in \Gamma \}$
from $\psi_2 | \Gamma$ for each $\Gamma \in \hat{\xi}_2$, and so on.  That
is, we define:
\begin{eqnarray*}
 \psi_i(z) &=& \frac{\tilde{\psi}_{i-1}(z)}{J\hat{F}(\hat{F}^{i-1}z)}; \\
 \varepsilon_{i,z} &=& \varepsilon_i \inf_{w \in \hat{\xi}_i(z)}
   \psi_i(w) ;\\
 \tilde{\psi}_i(z) &=& \psi_i(z) - \varepsilon_{i,z} .
\end{eqnarray*}
We show that under certain conditions on $(\varepsilon_i)$, the sequence
$\{ \tilde{\psi}_i \}$ satisfies a uniform bound on the ratios of its
values for nearby points.

A similar proposition to the following was obtained simultaneously but
independently by Holland (\cite{H1}); there, the emphasis was on the
regularity of the Jacobian.

\begin{prop}
\label{Ylemma3}
Suppose $(\varepsilon_i') \leq \frac{1}{2}$
is a sequence with the property that
\begin{eqnarray}
 \label{E1}
 v_i(\Phi) \prod_{j=1}^i (1+ \varepsilon_j') \leq K_0 ,
\end{eqnarray}
\begin{eqnarray}
 \label{E2}
 \sum_{j=1}^i \left( \prod_{k=j}^i (1 + \varepsilon_k') \right)
   \beta^{i-j+1} \leq K_0
\end{eqnarray}
are both satisfied for some sufficiently large constant $K_0$ allowed to
depend only on $F$ and $v_0(\Phi)$.  Then there exist
$\bar{\delta}<1$
and $\bar{C}>0$ each depending only on $F$, $C_\Phi$ and $v_0(\Phi)$
such that if we choose $\varepsilon_i = \bar{\delta}
\varepsilon_i'$ for each $i$, then for all $w,z$ with $\hat{s}(w,z) \geq
i \geq 1$,
\[ \left| \log \frac{\tilde{\psi}_i(w)}{\tilde{\psi}_i(z)} \right| \leq
   \bar{C} .\]

\end{prop}

\emph{Proof:} Suppose we are given such a sequence $(\varepsilon_i')$
and assume that for each $i$ we have
\begin{eqnarray}
\label{propassumpt}
 \left| \log \frac{\tilde{\psi}_i(w)}{\tilde{\psi}_i(z)} \right| \leq
   (1 + \varepsilon_i') \left| \log \frac{\psi_i(w)}{\psi_i(z)} \right|
\end{eqnarray}
for every $w,z$ with $\hat{s}(w,z) \geq i$.  We shall see that we may achieve
this by a suitable choice of $(\varepsilon_i)$.  We note that when
$\hat{s}(w,z) \geq i$,
\[ \left| \log \frac{\psi_i(w)}{\psi_i(z)} \right| \leq \left| \log
   \frac{\tilde{\psi}_{i-1}(w)}{\tilde{\psi}_{i-1}(z)} \right| + C_{\hat{F}}
   \beta^{\hat{s}(w,z)-i} ,\]
so
\[ \left| \log \frac{\tilde{\psi}_i(w)}{\tilde{\psi}_i(z)} \right| \leq
   (1 + \varepsilon_i') \left\{ \left| \log \frac{\tilde{\psi}_{i-1}(w)}
   {\tilde{\psi}_{i-1}(z)} \right| + C_{\hat{F}}
   \beta^{\hat{s}(w,z)-i} \right\} . \]
Applying this inductively, we obtain the estimate
\begin{eqnarray*} \left| \log \frac{\tilde{\psi}_i(w)}{\tilde{\psi}_i(z)}
   \right| &\leq& C_\Phi \left( \prod_{j=1}^i (1 + \varepsilon_j') \right)
   v_{\hat{s}(w,z)} (\Phi) \\ 
 &+& (1 + \varepsilon_i')C_{\hat{F}} \beta^{\hat{s}(w,z)-i} +
   (1 + \varepsilon_i')(1 + \varepsilon_{i-1}')C_{\hat{F}}
   \beta^{\hat{s}(w,z)-(i-1)} \\
 &+&\ldots + \prod_{j=1}^i (1 + \varepsilon_j') C_{\hat{F}}
   \beta^{\hat{s}(w,z)-1} . \end{eqnarray*}
We see this is bounded above by the constant
$ (C_\Phi + \beta^{-1} C_{\hat{F}}) K_0 =: \bar{C}$.  So we have
\[ \left| \log \frac{\psi_i(w)}{\psi_i(z)} \right| \leq \bar{C} +
   C_{\hat{F}} \textrm{ for } \hat{s}(w,z) \geq i .\]
It remains to examine the choice of sequence
$(\varepsilon_i)$ necessary for (\ref{propassumpt})
to hold.  For now,
let $(\varepsilon_i)$ be some sequence with $\varepsilon_i \leq
\varepsilon_i'$ for each $i$.  Let $\Gamma = \hat{\xi}_i(w) =
\hat{\xi}_i(z)$ and write $\varepsilon_{i,\Gamma} := \varepsilon_{i,w} =
\varepsilon_{i,z}$.  Then

\begin{eqnarray*}
  \left| \log \frac{\tilde{\psi}_i(w)}{\tilde{\psi}_i(z)} \right| - \left|
   \log \frac{\psi_i(w)}{\psi_i(z)} \right| 
 &\leq& \left| \log \left(\frac{\psi_i(w)-\varepsilon_{i,\Gamma}}{\psi_i(z)-
   \varepsilon_{i,\Gamma}} \cdot \frac{\psi_i(z)}{\psi_i(w)} \right)
   \right| \\
 &=& \left| \log \left( 1 + \frac{\varepsilon_{i,\Gamma} \psi_i(w) -
   \varepsilon_{i,\Gamma} \psi_i(z)}{(\psi_i(z) - \varepsilon_{i,\Gamma})
   \psi_i(w)} \right) \right| \\
 &=& \left| \log \left( 1 + \frac{\frac{\varepsilon_{i,\Gamma}}{\psi_i(z)} -
   \frac{\varepsilon_{i,\Gamma}}{\psi_i(w)}}{1 -
   \frac{\varepsilon_{i,\Gamma}}{\psi_i(z)}} \right) \right| \\
 &\leq& C_1 \left| \frac{\frac{\varepsilon_{i,\Gamma}}{\psi_i(z)} -
   \frac{\varepsilon_{i,\Gamma}}{\psi_i(w)}}{1 -
   \frac{\varepsilon_{i,\Gamma}}{\psi_i(z)}} \right| \end{eqnarray*}
for some $C_1>0$.

We see that $0 \leq \frac{\varepsilon_{i,\Gamma}}{\psi_i(w)} \leq
\varepsilon_i$ for all $w \in \Gamma$, and

\[ \frac{\frac{\varepsilon_{i,\Gamma}}{\psi_i(z)} -
   \frac{\varepsilon_{i,\Gamma}}{\psi_i(w)}}{1 -
   \frac{\varepsilon_{i,\Gamma}}{\psi_i(z)}} \geq - \varepsilon_i >
   - \frac{1}{2} , \]
so $C_1$ may be chosen so as not to depend on anything.  Continuing from the
estimate above,
\begin{eqnarray*}
 \left| \log \frac{\tilde{\psi}_i(w)}{\tilde{\psi}_i(z)} \right| -
   \left| \log \frac{\psi_i(w)}{\psi_i(z)} \right| &\leq& C_1
   \frac{\varepsilon_{i,\Gamma}}{\psi_i(z)} \left| 1 -
   \frac{\psi_i(z)}{\psi_i(w)} \right|
    \frac{1}{1 - \frac{\varepsilon_{i,\Gamma}}{\psi_i(z)}} \\
 &\leq& C_1 C_2 \frac{\varepsilon_i}{1- \varepsilon_i} \left| \log
   \frac{\psi_i(w)}{\psi_i(z)} \right| ,\end{eqnarray*}
where $C_2$ may be chosen independently of $i,w,z$ since $\frac{\psi_i(w)}
{\psi_i(z)} \geq e^{-(\bar{C} + C_{\hat{F}})}$, provided that at each stage we
choose $\varepsilon_i$ small enough that
$C_1 C_2 \frac{\varepsilon_i}{1 - \varepsilon_i} \leq \varepsilon_i'$.

We confirm that it is sufficient to take $\varepsilon_i = \bar{\delta}
\varepsilon_i'$ for small enough $\bar{\delta}>0$.  This means
\[ C_1 C_2 \frac{\varepsilon_i}{1-\varepsilon_i}= C_1 C_2 \frac{\bar{\delta}
   \varepsilon_i'}{1-\bar{\delta} \varepsilon_i'} < \frac{ C_1 C_2
   \bar{\delta}} {1 - \bar{\delta}} \varepsilon_i ' ,\]
so taking $\bar{\delta}= \frac{1}{1+C_1 C_2}$ is sufficient.
\hfill $\square$

\section{Choosing a sequence $(\varepsilon_i)$}

Having shown that it is sufficient for our purposes for the sequence
$(\varepsilon_i')$ to satisfy conditions (\ref{E1}) and (\ref{E2}), we now
consider how we might choose a sequence $(\varepsilon_i)$ which, subject to
these conditions, decreases as slowly as possible.  Having chosen a
sequence, we shall then estimate the rate of convergence this gives us.

\begin{lemma}
\label{choicelem}
Given a sequence $v_i(\Phi)$, there exists a sequence
$(\varepsilon_i') \leq \frac{1}{2}$ satisfying (\ref{E1}) and (\ref{E2})
such that for
$\varepsilon_i = \bar{\delta} \varepsilon_i'$, and any $K>1$,
\[ \prod_{j=1}^i \left( 1 - \frac{\varepsilon_j}{K} \right) \leq C \max
   \left( v_i(\Phi)^{\frac{\bar{\delta}}{K}}, \theta^i \right) \]
for some $\theta<1$ depending only on $F$, and some $C>0$.
\end{lemma}

\emph{Proof:} We start by defining a sequence $(v_i^*)>0$ as follows:
we let $v_0^* = v_0(\Phi)$ and $v_i^* = \max
(v_i(\Phi),c v_{i-1}^*)$, where $c$ is some constant such that
$ \exp \{ - \min (\frac{1}{2},\beta^{-1}-1) \} < c < 1$.
We claim that $v_i^* = \mathcal{O} (v_i(\Phi))$ unless
$v_i(\Phi)$ decays exponentially fast, in which case $v_i^*$ decays
at some (possibly slower) exponential rate.  To see this, suppose otherwise,
in the case where $v_i(\Phi)$ decays slower than any exponential
speed.  Then for large $i$ certainly $v_i^* > v_i(\Phi)$, and so
$v_i^* = cv_{i-1}^*$ for large $i$, and $(v_i^*)$ decays exponentially fast.
But this means $v_i(\Phi)$ decays exponentially fast, which is a
contradiction.

Let us now choose $\varepsilon_i' = \log \frac{v_{i-1}^*}{v_i^*}$.
(We ignore the trivial case $v_0(\Phi)=0$.)
We see
that all terms are small enough that (\ref{E2}) is satisfied, and in
particular, $\varepsilon_i \leq \frac{1}{2}$.
Furthermore, 
\begin{eqnarray*}
   v_i(\Phi) \prod_{j=1}^i (1 + \varepsilon_j') &\leq& v_i^* \exp
   \left\{ \sum_{j=1}^i \varepsilon_j' \right\} \\
  &=& v_i^* \exp \{ \log v_0^* - \log v_i^* \} \\
  &=& v_0(\Phi) .\end{eqnarray*}
For any $K>1$,
\begin{eqnarray*}
   \prod_{j=1}^i \left( 1 - \frac{\varepsilon_j}{K} \right) &\leq&
      \exp \left\{ - \frac{\bar{\delta}}{K} \sum_{j=1}^i \varepsilon_j'
      \right\} \\
  &=& \exp \left\{ - \frac{\bar{\delta}}{K} (\log v_0^* - \log v_i^*)
   \right\} \\
  &=& (v_0(\Phi)^{-1} v_i^*)^{\frac{\bar{\delta}}{K}} . \end{eqnarray*}

If $(v_i^*)$ decays exponentially fast, we get an exponential bound.
Otherwise, $v_i^{* \frac{\bar{\delta}}{K}} = \mathcal{O}
(v_i(\Phi)^\frac{\bar{\delta}}{K})$.

\section{Convergence of measures}

\label{convmeas}

We introduce a
sequence of measure densities $\hat{\Phi}_0 \equiv \Phi \geq \hat{\Phi}_1
\geq \hat{\Phi}_2 \geq \ldots$ corresponding to the sequence $\{
\tilde{\psi}_i \}$ in the following way:
\[ \hat{\Phi}_i(z):= \tilde{\psi}_i(z) J\hat{F}^{i}(z) .\]

\begin{lemma}
Given a sequence $(\varepsilon_i)=(\bar{\delta} \varepsilon_i')$ satisfying
the assumptions of Proposition \ref{Ylemma3}, there exists $K>1$ dependent
only on $F$, $C_\Phi$ and $v_0(\Phi)$ such that for all $z \in
\Delta \times \Delta$, $i \geq 1$,
\[ \hat{\Phi}_i(z) \leq \left( 1 - \frac{\varepsilon_i}{K} \right)
   \hat{\Phi}_{i-1}(z) .\]
\end{lemma}

\emph{Proof:}  If we fix $i \geq 1$, $\Gamma \in \hat{\xi}_i$, and $w,z \in
\Gamma$, then by Proposition \ref{Ylemma3} we have
\[ \frac{\hat{\Phi}_i (w)}{J\hat{F}^i(w)} \leq \bar{C}_0
   \frac{\hat{\Phi}_i (z)} {J\hat{F}^i (z)} \]
where $\bar{C}_0 = e^{\bar{C}}>1$.  From Sublemma \ref{sbreg}, we have
\[  \frac{1}{J\hat{F}(\hat{F}^{i-1} w)} \leq e^{C_{\hat{F}}} \frac{1}
   {J\hat{F} (\hat{F}^{i-1} z)} ,\]
so
\[ \frac{ \hat{\Phi}_{i-1}(w)}{J\hat{F}^i(w)} = \frac{\hat{\Phi}_{i-1}(w)}
   {J\hat{F}^{i-1}(w)} \cdot \frac{1}{J\hat{F} (\hat{F}^{i-1} w)}
   \leq \bar{C}_0 e^{C_{\hat{F}}} \frac{\hat{\Phi}_{i-1}(z)}{J\hat{F}^i(z)}
   .\]

Now we obtain a relationship between $\hat{\Phi}_i$ and $\hat{\Phi}_{i-1}$ by
writing
\begin{eqnarray*}
   \hat{\Phi}_i (z) &=& (\psi_i(z) - \varepsilon_{i,z}) J\hat{F}^i (z) \\
 &=& \left( \frac{\tilde{\psi}_{i-1}(z)}{J\hat{F}(\hat{F}^{i-1}z)} -
   \varepsilon_i \inf_{w \in \hat{\xi}_i(z)} \frac{\tilde{\psi}_{i-1}(w)}
   {J\hat{F} (\hat{F}^{i-1} w)} \right) J \hat{F}^i (z) \\
 &=& \left( \frac{\hat{\Phi}_{i-1}(z)}{J\hat{F}^i(z)} - \varepsilon_i
   \inf_{w \in \hat{\xi}_i (z)} \frac{\hat{\Phi}_{i-1}(w)}{J\hat{F}^i (w)}
   \right) J \hat{F}^i (z) . \end{eqnarray*}

So for any $z \in \Delta \times \Delta$ we have
that
\begin{eqnarray*}
   \hat{\Phi}_i(z) &\leq& \left( \frac{\hat{\Phi}_{i-1}(z)}{J\hat{F}^i (z)} -
   \frac{\varepsilon_i}{K} \frac{\hat{\Phi}_{i-1}(z)}{J\hat{F}^i (z)} \right)
   J \hat{F}^i (z) \\
 &=& \left( 1 - \frac{\varepsilon_i}{K} \right) \hat{\Phi}_{i-1}(z) 
   , \end{eqnarray*}
where $K=\bar{C}_0 e^{C_{\hat{F}}}$.
  \hfill $\square$

The above lemma gives an estimate on the total mass of $\hat{\Phi}_i$ for
each $i$.  To obtain an estimate for the difference between $F^n_* \lambda$
and $F^n_* \lambda'$, we must use this, and also take into account the
length of the simultaneous return time $T$.

\begin{lemma}
\label{Ylemma4}
For all $n>0$,
\[ \left| F^n_* \lambda - F^n_* \lambda' \right| \leq 2 P \{ T>n \} + 2
   \sum_{i=1}^n
   \left( \prod_{j=1}^i \left( 1 - \frac{\varepsilon_j}{K} \right) \right)
   P \{ T_i \leq n < T_{i+1} \} ,\]
where $K$ is as in the previous lemma.

\end{lemma}

\emph{Proof:}
We define a sequence $\{ \Phi_i \}$ of measure densities, corresponding to
the measure unmatched at time $i$ with respect to $F \times F$.  We shall
often write $\Phi_i (m \times m)$, say, to refer to the measure which has
density $\Phi_i$ with respect to $m \times m$.
For $z \in \Delta \times \Delta$ we let $\Phi_n(z) =
\hat{\Phi}_i(z)$, where $i$ is the largest integer such that $T_i(z) \leq n$.
Writing $\Phi = \Phi_n + \sum_{k=1}^n ( \Phi_{k-1} - \Phi_k )$, we have
\begin{eqnarray}
 \left| F^n_* \lambda - F^n_* \lambda' \right| &=& \left| \pi_*
   (F \times F)^n_* ( \Phi (m \times m)) - \pi_*' (F \times F)^n_* (\Phi
   (m \times m)) \right| \nonumber \\
 &\leq& \left| \pi_* (F \times F)^n_* (\Phi_n (m \times m)) - \pi_*' (F \times
   F)^n_* (\Phi_n (m \times m)) \right| \nonumber \\
\label{coorddiffeqn}
   &&+ \sum_{k=1}^n \left| (\pi_* - \pi_*') \left[ (F \times F)^n_* ((\Phi_{k-1}
   - \Phi_k)(m \times m)) \right] \right| . \end{eqnarray}

The first term is clearly $\leq 2 \int \Phi_n d(m \times m)$.  
Our construction should ensure the remaining terms are zero, since we have
arranged that the measure we subtract is symmetric in the two
coordinates.  To confirm this, we
partition $\Delta \times \Delta$ into regions on which each $T_m$ is
constant, at least while $T_m <n$.

Consider the family of sets $A_{k,i}$, $i,k
\in \mathbb{N}$, where $A_{k,i}:= \{ z \in \Delta \times \Delta: T_i(z)=k
\}$.  Clearly, each $A_{k,i}$ is a union of elements of $\hat{\xi}_i$, and
for any fixed $k$ the sets $A_{k,i}$ are pairwise disjoint.  It is also clear
that on any $A_{k,i}$, $\Phi_{k-1} - \Phi_k \equiv \hat{\Phi}_{i-1} -
\hat{\Phi}_i$, and for any $k$, $\Phi_{k-1} \equiv \Phi_k$ on $\Delta \times
\Delta - \cup_i A_{k,i}$.  So for each $k$,
\[ \pi_* (F \times F)^n_* ((\Phi_{k-1} - \Phi_k) (m \times m)) \]
\[ = \sum_i \sum_{\Gamma \in (\hat{\xi}_i|A_{k,i})} F^{n-k}_* \pi_* (F \times F)^{T_i}_*
   ((\hat{\Phi}_{i-1} - \hat{\Phi}_i) ((m \times m) | \Gamma)) .\]

We show that this measure is unchanged if we replace $\pi$ with $\pi'$ in the
last expression.  Let $E \subset \Delta$ be an arbitrary
measurable set, and fix some $\Gamma \in \hat{\xi}_i | A_{k,i}$.
  Then
\[ \pi_* \hat{F}^i_* ((\hat{\Phi}_{i-1} - \hat{\Phi}_i) ((m \times m)|
   \Gamma))(E) \]
\[ = \hat{F}^i_* \left( \left( \varepsilon_i J\hat{F}^i \inf_{w \in
   \Gamma} \frac{\hat{\Phi}_{i-1}(w)}{J\hat{F}^i(w)} 
   \right) \left( (m \times m) | \Gamma \right) \right) (E \times \Delta) \]
\[ = \left( \varepsilon_i C J\hat{F}^i (m \times m) \right)
   (\hat{F}^{-i} (E \times \Delta) \cap \Gamma) \]
where $C$ is constant on $\Gamma$.   This equals
\[ \int_{\hat{F}^{-i} (E \times \Delta) \cap \Gamma} \varepsilon_i C
   J\hat{F}^i d(m \times m) = \varepsilon_i C (m \times m) (E \times
   \Delta) .\]
Since $(m \times m) (E \times \Delta) = (m \times m) (\Delta \times E)$, the
terms of the sum in (\ref{coorddiffeqn}) all have zero value, as claimed.

Now \[ \int \Phi_n d(m \times m) = \sum_{i=0}^\infty
\int_{\{T_i \leq n < T_{i+1}\}} \Phi_n d(m \times m); \]
in fact, since $T_i \geq i$, all terms of the series are zero for $i>n$.
For $1 \leq i \leq n$,
\[ \int_{ \{ T_i \leq n < T_{i+1} \}} \Phi_n = \int_{ \{ T_i \leq n <
   T_{i+1} \}} \hat{\Phi}_i
   \leq \int_{ \{ T_i \leq n < T_{i+1} \} } \prod_{j=1}^i \left( 1 -
   \frac{\varepsilon_j}{K} \right) \Phi .\]

The estimate claimed for $|F^n_* \lambda - F^n_* \lambda'|$
follows easily. \hfill $\square$

Finally we state a simple relationship between $P \{ T>n \}$ and
$(m \times m) \{ T>n \}$.  From now on we shall use the convention that
$P \{ \textrm{condition} | \Gamma \} :=
\frac{1}{P(\Gamma)} P \{ x \in \Gamma: x \textrm{ satisfies
condition} \}$.

\begin{sublemma}
\label{Ysb4}
There exists $\bar{K}>0$ depending only on $C_\Phi$ and
$v_0(\Phi)$ s.t. $\forall i \geq 1$, $\forall \Gamma \in \hat{\xi}_i$,
\[ P \{ T_{i+1}-T_i >n | \Gamma \} \leq \bar{K} (m_0 \times m_0) \{ T>n \}
   .\]
The dependence of $\bar{K}$ on $P$ may be removed entirely if we take only
$i \geq$ some $i_0(P)$.
\end{sublemma}

\emph{Proof:} Let $\mu = \frac{1}{P(\Gamma)} \hat{F}^i_* (P| \Gamma)$.  We
see that $P \{ T_{i+1} - T_i >n | \Gamma \} = \mu \{ T>n \}$.  We prove a
distortion estimate for $\frac{d\mu}{d(m_0 \times m_0)}$, using the estimates
of Sublemma \ref{sbprodreg}.  Let $w,z \in \Delta_0 \times \Delta_0$ and let
$w_0 ,z_0 \in \Gamma$ be such that $\hat{F}^i w_0 = w$, $\hat{F}^i z_0 = z$.
Then
\begin{eqnarray*}
 \left| \log \frac{\frac{d\mu}{d(m \times m)} (w)}{\frac{d\mu}{d(m \times 
   m)}(z)} \right| &\leq& \left| \log \frac{\Phi(w_0)}{\Phi(z_0)} \right| +
   \left| \log \frac{J\hat{F}^i w_0}{J\hat{F}^i z_0} \right| \\
  &\leq& C_\Phi v_i(\Phi) + C_{\hat{F}} . \end{eqnarray*}

This gives
\[ \frac{d\mu}{d(m \times m)} \leq \frac{e^{(C_\Phi v_i(\Phi) +
   C_{\hat{F}})}}{(m_0 \times m_0) (\Delta_0 \times \Delta_0)} \]
and hence the result follows.

\section{Combinatorial estimates}
\label{combinatorics}

In Lemma \ref{Ylemma4} we have given the main estimate involving
$P$, $T$, and the sequence $(\varepsilon_i)$.  It
remains to relate $P$ and $T$ to the sequence $m_0 \{
R>n \}$.
Primarily, this involves estimates relating the sequences $P \{ T>n \}$
and $m_0 \{ R>n \}$.  We shall state only some key estimates of the proof,
referring the reader to \cite{Y2} for full details.  Our statements differ
slightly, as the estimates of \cite{Y2}  are stated in terms of $m \{
\hat{R}>n \}$;
they are easily reconciled by noting that $m \{ \hat{R}>n \} = \sum_{i>n}
m_0 \{ R>i \}$.
(As earlier, $\hat{R} \geq 0$ is the first arrival time to $\Delta_0$.)

\begin{prop}
\label{Ycomb}
\begin{enumerate}
\item If $m_0 \{ R>n \} = \mathcal{O} (\theta^n)$ for some $0<\theta<1$,
   then $P \{ T>n \} = \mathcal{O} (\theta_1^n)$, for some $0<\theta_1<1$.
   Also,  for sufficiently small $\delta_1>0$,
\[ P \{ T_i \leq n < T_{i+1} \}
   \leq C \theta'^n \textrm{ for } i \leq \delta_1 n , \]
for some $0<\theta'<1$, $C>0$ independent of $i$.  The constants $\theta_1$,
$\theta'$, $\delta_1$ may all be chosen independently of $P$.
\item If $m_0 \{ R>n \} = \mathcal{O} (n^{-\alpha})$ for some $\alpha>1$,
   then $P \{ T>n \} = \mathcal{O}(n^{1-\alpha})$.
\end{enumerate}
\end{prop}

This proposition follows from estimates involving the combinatorics of the
intermediate stopping times $\{ \tau_i \}$.
Let us make explicit a key sublemma used in the proofs, concerning the
regularity of
the pushed-forward measure densities $\frac{dF^n_* \lambda}{dm}$; for the
rest of the argument we refer to \cite{Y2}, as the changes are minor.

\begin{sublemma}
\label{sbreg}
For any $k>0$, let $\Omega \in \bigvee_{i=0}^{k-1} F^{-i}\eta$ be s.t.
$F^k \Omega = \Delta_0$.  Let $\mu = F^k_* (\lambda | \Omega)$.  Then
$\forall x,y \in \Delta_0$, we have
\[ \left| \frac{\frac{d\mu}{dm}(x)}{\frac{d\mu}{dm}(y)} -1 \right| \leq C_0 \]
for some $C_0(\lambda)$, where the dependence on $\lambda$ is only on
$v_0(\lambda)$ and $C_\lambda$, and may be removed entirely if
we only consider $k \geq$ some $k_0(\lambda)$.
\end{sublemma}

\emph{Proof:} Let $\varphi = \frac{d\lambda}{dm}$, fix $x,y \in \Delta_0$,
and let $x_0,y_0$ be the unique points in $\Omega$ such that
$F^k x_0 =x$, $F^k y_0 = y$.  We note that $\frac{d\mu}{dm}(x) = \varphi(x_0)
\cdot \frac{dF^k_* (m|\Omega)}{dm}(x_0)$.  
So
\begin{eqnarray*}
  &&\left| \frac{\varphi (x_0)}{JF^k x_0} \cdot \frac{JF^k y_0}{\varphi (y_0)}
   -1 \right| = \frac{JF^k y_0}{\varphi (y_0)} \left| \frac{\varphi (x_0)}
   {JF^k x_0} - \frac{\varphi (y_0)}{JF^k y_0} \right| \\
  &\leq & \frac{JF^k y_0}{\varphi (y_0)} \left\{ \varphi (x_0) \left|
   \frac{1}{JF^k x_0} - \frac{1}{JF^k y_0} \right| + \frac{1}{JF^k y_0}
   | \varphi (x_0) - \varphi (y_0) | \right\} \\
 &\leq& \frac{\varphi (x_0)}{\varphi(y_0)} \left| \frac{JF^k y_0}{JF^k x_0}-1
   \right| + \left| \frac{\varphi (x_0)}{\varphi (y_0)} -1 \right| \\
 &\leq& (1 + C_\varphi v_j(\varphi)) C' + C_\varphi
   v_j(\varphi) \leq (1 + C_\varphi v_0(\varphi))C' +
   C_\varphi v_0(\varphi) ,\end{eqnarray*}
where $j$ is the number of visits to $\Delta_0$ up to time $k$.  Clearly the
penultimate bound can be made independent of $\lambda$ for
$j \geq$ some $j_0 (\lambda)$.

The following result combines the estimates above with those of the previous
section.

\begin{prop}
\label{myprop}
\begin{enumerate}
\item When $m_0 \{ R>n \} \leq C_1 \theta^n$,
\[ \left| F^n_* \lambda - F^n_* \lambda' \right| \leq C \theta'^n + 2
   \sum_{i=[\delta_1 n]+1}^n \left(
   \prod_{j=1}^i (1 - \frac{\varepsilon_j}{K}) \right) \]
for some $0<\theta'<1$ and sufficiently small $\delta_1$.
\item When $m_0 \{ R>n \} \leq C_1 n^{-\alpha}$, $\alpha>1$,
\begin{eqnarray*}
 \left| F^n_* \lambda - F^n_* \lambda' \right| \leq 2 C_1 n^{1-\alpha} &+&
   C n^{1-\alpha}
   \sum_{i=1}^{[\delta_1 n]} i^{\alpha} \prod_{j=1}^i
   \left( 1 - \frac{\varepsilon_j}{K} \right) \\
   &+& \sum_{i=[\delta_1 n]+1}^n
   \prod_{j=1}^i (1 - \frac{\varepsilon_j}{K})  \end{eqnarray*}
for sufficiently small $\delta_1$.
\end{enumerate}
\end{prop}

\emph{Proof:}
In the first case, Proposition \ref{Ycomb} and Lemma \ref{Ylemma4} tell us
that
\[ \left| F^n_* \lambda - F^n_* \lambda' \right| \leq C \theta_0^n +
   2 \sum_{i=1}^{[\delta_1 n]} P \left\{ T_i \leq n < T_{i+1} \right\} +
   2 \sum_{i=[\delta_1 n]+1}^n \prod_{j=1}^i \left( 1 -
   \frac{\varepsilon_j}{K} \right) \]
for any $0<\delta_1<1$, for some $0<\theta_0<1$.  For sufficiently small
$\delta_1$, the middle term is $\leq C [\delta_1 n] \theta'^n$, which decays
at some exponential speed in $n$.

In the second case, for any $0<\delta_1<1$ we have
\begin{eqnarray*}
 \left| F^n_* \lambda - F^n_* \lambda' \right| \leq C n^{1-\alpha} &+&
   2 \sum_{i=1}^{[\delta_1 n]} \left( \prod_{j=1}^i \left( 1 -
   \frac{\varepsilon_j}{K} \right) \right) P \left\{ T_i \leq n < T_{i+1}
   \right\} \\ &+& \sum_{i=[\delta_1 n] + 1}^n \prod_{j=1}^i \left( 1 -
   \frac{\varepsilon_j}{K} \right) . \end{eqnarray*}

We estimate the middle term by noting that
\begin{eqnarray*}
 P \{ T_i \leq n < T_{i+1} \} &\leq& \sum_{j=0}^i P \left\{ T_{j+1} - T_j >
   \frac{n}{i+1} \right\} \\
 &\leq& \bar{K} (i+1) (m \times m) \left\{ T> \frac{n}{i+1} \right\} \\
 &\leq& C n^{1 - \alpha} (i+1)^\alpha \textrm{ for some } C>0.
 \end{eqnarray*}
For the last step, note that Proposition \ref{Ycomb} applies to the
normalisation of $(m \times m)$ to a probability measure.

\section{Specific regularity classes}

We now combine all of our intermediate estimates to obtain a rate of
decay of correlations in the specific cases mentioned in Theorem \ref{Ythm2}.
First, we set $\zeta = \frac{\bar{\delta}}{K}$, which can be seen to
depend only on $F$, $C_\Phi$ and $v_0(\Phi)$.  Throughout this
section, we shall let $C$ denote a generic constant, allowed to depend only
on $F$ and $\Phi$, which may vary between expressions.

\subsection{Exponential return times}
In this subsection, we suppose that $m_0 \{ R>n \} = \mathcal{O} (\theta^n)$, and
hence $m \{ \hat{R}>n \} = \mathcal{O} (\theta^n)$.

\emph{Class (V1):} Suppose $v_i(\Phi) = \mathcal{O}(\theta_1^i)$
for some $\theta_1<1$.  By Lemma \ref{choicelem} we may take
$(\varepsilon_i)$ such that
conditions (\ref{E1}) and (\ref{E2}) are satisfied, and

\[ \prod_{j=1}^i \left( 1 - \frac{\varepsilon_j}{K} \right) = \mathcal{O}
  ( \theta_2^i ) \]
for some $0<\theta_2<1$.
Applying Proposition \ref{myprop} we have
\[ \left| F^n_* \lambda - F^n_* \lambda' \right| \leq C \theta'^n + C
   \sum_{i> [\delta_1 n]} \theta_2^i \]
for some $\theta'<1$ and sufficiently small $\delta_1>0$.  This give the
required exponential bound in $n$.

\emph{Class (V2):} Suppose $v_i(\Phi ) = \mathcal{O} \left(
e^{-i^\gamma} \right)$, for some $\gamma \in (0,1)$.  Then there exists
$(\varepsilon_j)$ such that
\[ \prod_{j=1}^i \left( 1 - \frac{\varepsilon_j}{K} \right) = \mathcal{O}
   (e^{-\zeta i^\gamma}) .\]
So 
\[ \left| F^n_* \lambda - F^n_* \lambda' \right| \leq C \theta'^n + C
   \sum_{i>[\delta_1 n]} e^{-\zeta i^\gamma} .\]

We see that $e^{-\zeta i^\gamma} = \mathcal{O} (e^{-i^{\gamma'}})$ for every
$0<\gamma'<\gamma$, and it is well known that the sum is of order
$e^{-n^{\gamma''}}$ for every $0<\gamma''<\gamma'$.

\emph{Class (V3):} Suppose $v_i(\Phi)=\mathcal{O} \left(e^{-(\log
i)^\gamma} \right)$ for some $\gamma>1$.  We may take $(\varepsilon_j)$
such that
\[ \prod_{j=1}^i \left(1 - \frac{\varepsilon_j}{K} \right) = \mathcal{O}
   (e^{-\zeta  (\log i)^\gamma}) .\]
So
\[ \left| F^n_* \lambda - F^n_* \lambda' \right| \leq C \theta'^n +
   C \sum_{i>[\delta_1 n]} e^{-\zeta (\log i)^\gamma} . \]
It is easy to show that $e^{-\zeta (\log i)^\gamma} = 
\mathcal{O} \left(e^{-(\log i)^{\gamma'}} \right)$
for every $0<\gamma'<\gamma$.  So the sum is of order
$\mathcal{O} (e^{-(\log n)^{\gamma''}})$ for every $0<\gamma''<\gamma$.

\emph{Class (V4):} Suppose $v_i(\Phi) = \mathcal{O}(i^{-\gamma})$
for some $\gamma>\frac{1}{\zeta}$.  Then we can take $(\varepsilon_j)$
such that
\[ \prod_{j=1}^i \left( 1 - \frac{\varepsilon_j}{K} \right) \leq C
   i^{-\zeta \gamma} . \]
So
\[ \left| F^n_* \lambda - F^n_* \lambda' \right| \leq C \theta'^n + C
   \sum_{i=[\delta_1 n +1]}^n i^{-\zeta \gamma} = \mathcal{O} (n^{1-\zeta
   \gamma}) .\]

\subsection{Polynomial return times}
Here we suppose $m_0 \{ R>n \} = \mathcal{O} (n^{-\alpha})$ for some
$\alpha>1$.
Suppose $v_n(\Phi) = \mathcal{O}
(n^{-\gamma})$, for some $\gamma > \frac{2}{\zeta}$.  We can take $(\varepsilon_i)$ such that
\[ \prod_{j=1}^i \left( 1 - \frac{\varepsilon_j}{K} \right) \leq C i^{-\zeta
\gamma} .\]
By Proposition \ref{myprop}, for some $\delta_1$,
\[ \left| F^n_* \lambda - F^n_* \lambda' \right| \leq 2 C n^{1-\alpha}
   + C n^{1-\alpha} \sum_{i=1}^{[\delta_1 n]} i^{\alpha -
   \zeta \gamma} + C \sum_{i=[\delta_1 n]+1}^\infty i^{-\zeta \gamma} .\]

The third term here is of order $n^{1-\zeta \gamma}$.  To estimate the second term, we consider
three cases.

\emph{Case 1: $\gamma> \frac{\alpha+1}{\zeta}$.}  Here, $\alpha - \zeta
\gamma < -1$, so the sum is bounded above independently of $n$, and the
whole term is $\mathcal{O}(n^{1-\alpha})$.

\emph{Case 2: $\gamma = \frac{\alpha+1}{\zeta}$.}  The sum is
\[ \sum_{i \leq [\delta_1 n]} i^{-1} \leq 1 + \int_1^{[\delta_1 n]} x^{-1}
   dx = 1 + \log [\delta_1 n] = \mathcal{O} (\log n) .\]
So the whole term is $\mathcal{O} (n^{1-\alpha} \log n)$.

\emph{Case 3: $\frac{2}{\zeta}< \gamma < \frac{\alpha+1}{\zeta}$.}  The sum
is of order $n^{\alpha+1 - \zeta \gamma}$, and so the whole term is
$\mathcal{O}(n^{2-\zeta \gamma})$.

\section{Decay of correlations}

\label{correlations}

Finally, we show how estimates for decay of correlations may be
derived directly from those for the rates of convergence of measures.

Let $\varphi,\psi \in \mathcal{L}^\infty (\Delta,m)$, as
in the statement of Theorem \ref{Ythm2}.  We write $\tilde{\psi}:= b (\psi +
a)$, where $a = 1 - \inf \psi$ and $b$ is such that $\int \tilde{\psi} d\nu
=1$.  We notice that $b \in \left[ \frac{1}{1+v_0(\psi)},1 \right]$,
and that $\inf \tilde{\psi}=b$, $\sup \tilde{\psi} \leq
1 + v_0(\psi)$.

Now let $\rho = \frac{d\nu}{dm}$, and let
$\lambda$ be the measure on $\Delta$ with $\frac{d\lambda}{dm} =
\tilde{\psi} \rho$.  We have
\begin{eqnarray*}
 \left| \int \left( \varphi \circ F^n \right) \psi d\nu - \int \varphi
   d\nu \int \psi d\nu \right|
 &=& \frac{1}{b} \left| \int \left( \varphi \circ F^n \right) \tilde{\psi}
   d\nu - \int \varphi d\nu \int \tilde{\psi} d\nu \right| \\
 &=& \frac{1}{b} \left| \int \varphi \, d \left( F^n_* (\tilde{\psi} \rho m)
   \right) - \int \varphi \rho dm \right| \\
 &\leq& \frac{1}{b} \int | \varphi | \left| \frac{dF^n_* \lambda}{dm} - \rho
   \right| dm \\
 &\leq& \frac{1}{b} \| \varphi \|_\infty \left| F^n_* \lambda - \nu \right|
   . \end{eqnarray*}

It remains to check the regularity of $\tilde{\psi} \rho$.  First,
\begin{eqnarray*}
 \left| \tilde{\psi}(x) \rho(x) - \tilde{\psi}(y) \rho(y) \right| &\leq&
   \left| \tilde{\psi}(x) \left( \rho(x) - \rho(y) \right) \right| +
   \left| \rho(y) \left( \tilde{\psi}(x) - \tilde{\psi}(y) \right) \right| \\
 &\leq& \| \tilde{\psi} \|_\infty \left| \rho(x) - \rho(y) \right| +
   \| \rho \|_\infty |b| \left| \psi(x) - \psi(y) \right| . \end{eqnarray*}

It can be shown that $\rho$ is bounded below by some positive constant, and
$v_n(\rho) \leq C \beta^n$.  (This is part of
the statement of Theorem 1 in \cite{Y2}.)  So $\tilde{\psi} \rho$ is bounded
away from zero, and
$v_n({\tilde{\psi}\rho}) \leq C \beta^n + C v_n(\psi)$,
where $C$ depends on $v_0(\psi)$.

Taking $\lambda'=\nu$, we see that $v_n(\Phi) \leq
C \beta^n + C v_n(\psi)$.
This shows that estimates for $| F^n_* \lambda -
F^n_* \lambda'|$ carry straight over to estimates for decay of correlations.
To check that the dependency of the constants is as we require, we note
that we can take
\[ C_\lambda = \frac{1}{\inf \frac{d\lambda}{dm}} =
   \frac{1}{\inf \tilde{\psi} \rho} 
   \leq \frac{1}{b \inf \rho} \leq \frac{1 + v_0(\psi)}{\inf
   \rho} . \]
So an upper bound for this constant is determined by $v_0(\psi)$.
Clearly $C_\Phi$ depends only on $F$ and an upper bound for $v_0(\psi)$,
and in particular these constants determine $\zeta =
   \frac{\bar{\delta}}{K}$.

\section{Central Limit Theorem}

We verify the Central Limit Theorem in each case for classes of
observables which give summable decay of autocorrelations (that is,
summable decay of correlations under the restriction $\varphi=\psi$).

A general theorem of Liverani (\cite{L2}) reduces in this context to the
following.

\begin{theorem}
Let $(X,\mathcal{F},\mu)$ be a probability space, and $T:X \circlearrowleft$
a (non-invertible) ergodic measure-preserving transformation.  Let $\varphi
\in \mathcal{L}^\infty (X,\mu)$ be such that $\int \varphi d\mu = 0$.  Assume
\begin{eqnarray}
 \label{clt1}
   &&\sum_{n=1}^{\infty} \left| \int (\varphi \circ T^n) \varphi d\mu \right|
   < \infty , \\
 \label{clt2}
   &&\sum_{n=1}^{\infty} ({\hat{T}}^{*n} \varphi)(x) \textrm{ is absolutely
   convergent for $\mu$-a.e. $x$} ,  \end{eqnarray}
where $\hat{T}^*$ is the dual of the operator $\hat{T}: \varphi \mapsto
\varphi \circ T$.
Then the Central Limit Theorem holds for $\varphi$
if and only if $\varphi$ is not a coboundary.
\end{theorem}

In the above, the dual operator $\hat{T}^*$ is the Perron-Frobenius
operator corresponding to $T$ and $\mu$, that is
\[ (\hat{T}^* \varphi) (x) = \sum_{y:Ty=x} \frac{\varphi(y)}{JT(y)} .\]
Of course the Jacobian $JT$ here is defined in terms of the measure $\mu$.

Let
$\varphi:\Delta \rightarrow \mathbb{R}$ be an observable which is not a
coboundary, and for which
\[ \mathcal{C}_n (\varphi,\varphi;\nu) = \left| \int (\varphi \circ F^n)
   \varphi d\nu - (\int \varphi d\nu)^2 \right| \]
is summable.  Let $\phi = \varphi
- \int \varphi d\nu$, so that $\int \phi d\nu = \int \varphi d\nu - \int
(\int \varphi d\nu) d\nu = 0$.  We shall show that $\phi$ satisfies the
assumptions of the theorem above.  It is straightforward to check that
$\mathcal{C}_n (\varphi,\varphi;\nu) = \mathcal{C}_n (\phi,\phi;\nu)
      = \left| \int (\phi \circ F^n) \phi d\nu
\right|$.  Hence condition (\ref{clt1}) above is satisfied for $\phi$.

Since $m$ and $\nu$ are equivalent measures, it suffices to verify the
condition in (\ref{clt2}) $m$-a.e.  The operator $\hat{F}^*$ is defined
in terms of the invariant measure, so for a measure $\lambda \ll m$ it
sends $\frac{d\lambda}{d\nu}$ to $\frac{dF_* \lambda}{d\nu}$.  By a change
of coordinates (or rather, of reference measure), we find that
\begin{eqnarray*} (\hat{F}^{*n} \phi) (x) = \frac{1}{\rho(x)}
   (\mathcal{P}^n (\phi \rho)) (x) ,\end{eqnarray*}
where $\mathcal{P}$ is the Perron-Frobenius operator with respect to $m$,
that is, the operator sending densities $\frac{d\lambda}{dm}$ to
$\frac{dF_* \lambda}{dm}$.

We shall now write $\phi$ as the difference of the densities of two
(positive) measures of similar regularity to $\phi$.  We let
$\tilde{\phi} = b(\phi +a)$, for some large $a$, with
$b>0$ chosen such that $\int \tilde{\phi} \rho dm =1$.  We define measures
$\lambda,\lambda'$ by
\[ \frac{d\lambda}{dm} = (b \phi + \tilde{\phi}) \rho , \, 
   \frac{d\lambda'}{dm} = \tilde{\phi} \rho .\]
It is straightforward to check this gives two probability measures, and
that 
\[ b^{-1} \left( \frac{d\lambda}{dm} - \frac{d\lambda'}{dm} \right) = \phi
   \rho . \]
As we showed in the previous section,
$v_n(\tilde{\phi} \rho) \leq C \beta^n + C v_n(\tilde{\phi})$ for
some $C>0$.  Also, $b\phi + \tilde{\phi} = b(2\phi+a)$, which is
bounded below by some positive constant, provided we choose sufficiently
large $a$.  We easily see $v_n(\lambda),v_n(\lambda') \leq C
v_n (\varphi)$.

We now follow the construction of the previous sections for these given
measures $\lambda$, $\lambda'$, and consider the sequence of densities
$\Phi_n$ defined in \S \ref{convmeas}.  We have
\[ F^n_* \lambda - F^n_* \lambda' = \pi_* (F \times F)^n_* (\Phi_n (m \times
   m)) - \pi_*' (F \times F)^n_* (\Phi_n (m \times m)) .\]
Let $\psi_n$ be the density of the first term with respect to $m$, and
$\psi_n'$ the density of the second.  Since $\mathcal{P}$ is a linear
operator, we see that 
\begin{eqnarray*} | \mathcal{P}^n (\phi \rho)| &=& b^{-1} \left| \frac{dF^n_*
   \lambda}{dm} - \frac{dF^n_* \lambda'}{dm} \right| \\
  &\leq& b^{-1} (\psi_n + \psi_n') \end{eqnarray*}
These densities have integral and distortion which are estimable by
the construction.  We know $\int \psi_n dm = \int \psi_n' dm = \int \Phi_n
d(m \times m)$.  In the cases we consider (sufficiently fast polynomial
variations) this is summable in $n$; notice that we have already used this
expression as a key upper bound for $\frac{1}{2} |F^n_* \lambda - F^n_*
\lambda'|$ (see Lemma \ref{Ylemma4}).
It remains to show that a similar condition holds pointwise, by showing that
$\psi_n$, $\psi_n'$ both have bounded distortion
on each $\Delta_l$, and hence $|F^n_* \lambda - F^n_* \lambda'|$ is
an upper bound for $\psi_n + \psi_n'$, up to some constant.  This follows
non-trivially from Proposition \ref{Ylemma3}, which gives a distortion bound
on $\{ \Phi_k \}$, and hence on $\{ \Phi_n \}$ when we restrict to elements
of a suitable partition.  The remainder of the argument is essentially no
different from that given in \cite{Y2}, and we omit it here.

\section{Applications}
\label{applications}

Having obtained estimates in the abstract framework of Young's tower,
we now discuss how these results may be applied to other settings.
First, we define formally what it means for a
system to \emph{admit} a tower.

Let $X$ be a finite dimensional compact Riemannian manifold, with
$\mathrm{Leb}$ denoting some Riemannian volume (\emph{Lebesgue measure}) on
$X$.
We say that a locally $C^1$ non-uniformly expanding system
$f:X \circlearrowleft$ \emph{admits a
tower} if there exists a subset $X_0 \subset X$, $\mathrm{Leb}(X_0)>0$, a
partition
(mod $\mathrm{Leb}$) $\mathcal{P}$ of $X_0$, and a return time function
$R:X_0
\rightarrow \mathbb{N}$ constant on each element of $\mathcal{P}$, such that

\begin{itemize}
\item for every $\omega \in \mathcal{P}$,
$f^R|\omega$ is an injection onto $X_0$;
\item $f^R$ and $(f^R |\omega)^{-1}$ are $\mathrm{Leb}$-measurable
functions, $\forall \omega \in \mathcal{P}$;
\item $\bigvee_{j=0}^\infty (f^R)^{-j} \mathcal{P}$ is the trivial
partition into points;
\item the volume derivative $\det Df^R$ is well-defined and non-singular
(i.e. $0<|\det Df^R|<\infty$) $\mathrm{Leb}$-a.e., and 
$\exists C>0$, $\beta<1$, such that $\forall
\omega \in \mathcal{P}$, $\forall x,y \in \omega$,
\[ \left| \frac{\det Df^R (x)}{\det Df^R (y)} -1 \right| \leq C
   \beta^{s(F^R x, F^Ry)} , \]
where $s$ is defined in terms of $f^R$ and $\mathcal{P}$ as before.
\end{itemize}

We say the system admits the tower $F:(\Delta,m) \circlearrowleft$ if the
base $\Delta_0 = X_0$, $m|\Delta_0 = \mathrm{Leb}|X_0$, and the tower is
determined
by $\Delta_0$, $R$, $F^R:=f^R$ and $\mathcal{P}$ as in \S \ref{towers}.  It
is easy to check that the usual assumptions of the tower hold, except
possibly for aperiodicity and finiteness.  In particular, $| \det Df^R|$
equals the Jacobian $JF^R$.

If $F:(\Delta,m) \circlearrowleft$ is a tower for $f$ as above, there exists
a projection $\pi:\Delta \rightarrow X$ we shall simply call the
\emph{tower projection},
which is a semi-conjugacy between $f$ and $F$; that is,
for $x \in \Delta_l$, with $x=F^l x_0$ for $x_0 \in
\Delta_0$, $\pi(x) := f^l (x_0)$.
In all the examples we have mentioned in \S \ref{examples}, the standard
tower constructions (as given in the papers we cited there)
provide us with a tower projection
$\pi$ which is H\"older-continuous with respect to 
the separation time $s$ on $\Delta$.  That is, given a Riemannian metric $d$,
in each case we have that
\begin{eqnarray}
\label{app0}
\exists \beta<1 \textrm{ such that for } x,y \in \Delta,
  \, d(\pi (x), \pi (y)) = \mathcal{O}(\beta^{s(x,y)}) .\end{eqnarray}
Note that the issue of the regularity of $\pi$ is often not mentioned
explicitly in the literature, but essentially follows from having good
distortion control for every iterate of the map.  (Formally, a tower
is only required to have good distortion for the return map $F^R$, which
is not sufficient.)

Given a system $f$ which admits a tower $F:(\Delta,m)
\circlearrowleft$ with projection $\pi$ satisfying (\ref{app0}), we show
how the observable classes $(R1-4)$ on $X$ correspond to the classes
$(V1-4)$ of observables on $\Delta$.  Recall that for given $\psi$,
\[ \mathcal{R}_\varepsilon (\psi) : = \sup \{ | \psi(x) - \psi(y) | :
   d(x,y) \leq \varepsilon \} . \]

Given a regularity for $\psi$ in terms of $\mathcal{R}_\varepsilon (\psi)$,
we estimate the regularity of $\psi \circ \pi$, which is an observable on
$\Delta$. 

\begin{lemma}
\begin{itemize}
\item If $\psi \in (R1,\gamma)$ for some $\gamma \in (0,1]$, then $\psi \circ \pi
   \in (V1)$;
\item if $\psi \in (R2,\gamma)$ for some $\gamma \in (0,1)$, then $\psi \circ \pi
   \in (V2,\gamma')$ for every $\gamma'<\gamma$;
\item if $\psi \in (R3,\gamma)$ for some $\gamma>1$, then $\psi \circ \pi \in
   (V3,\gamma')$ for every $\gamma'<\gamma$;
\item if $\psi \in (R4,\gamma)$ for some $\gamma>1$, then $\psi \circ \pi \in
   (V4,\gamma)$.
\end{itemize}
\end{lemma}

\emph{Proof:} The computations are entirely straightforward, so we shall just
make explicit the $(R4)$ case for the purposes of illustration.

Suppose $\mathcal{R}_\varepsilon (\psi) = \mathcal{O}
   (\left| \log \varepsilon \right|^{-\gamma})$, for some $\gamma>0$.
Then, taking $n$ large as necessary,
\begin{eqnarray*} v_n(\psi \circ \pi) &\leq& C | \log C \beta^n| ^{-\gamma}
   \\
   &=& C (n \log \beta^{-1} - \log C)^{-\gamma} \\
   &\leq& C (\frac{n}{2} \log \beta^{-1})^{-\gamma} \\
   &=& \mathcal{O} (n^{-\gamma}) .\end{eqnarray*}
\hfill $\square$

Let us point out that the condition (\ref{app0}) is not necessary for us
to apply these methods.  If we are given some weaker regularity on $\pi$,
the classes $(V1-4)$ shall simply correspond to some larger observable
classes on the manifold.
It remains to check that the semi-conjugacy $\pi$ preserves the
statistical properties we are interested in.

\begin{lemma}
Let $\nu$ be the mixing acip on $\Delta$ given by Theorem \ref{Ythm2}.
Given $\varphi,\psi:X \rightarrow \mathbb{R}$, let $\hat{\varphi} =
\varphi \circ \pi$, $\hat{\psi} = \psi \circ \pi$.  Then

\[ \mathcal{C}_n (\varphi,\psi;\pi_* \nu) = \mathcal{C}_n
   (\hat{\varphi},\hat{\psi};\nu) .\]

\end{lemma}

\begin{lemma}
Suppose the Central Limit Theorem holds for $(F,\nu)$ for some
observable $\varphi:\Delta \rightarrow \mathbb{R}$.
Then the Central Limit Theorem also holds for $(f,\pi_* \nu)$ for the observable
$\hat{\varphi} = \varphi \circ \pi$.
\end{lemma}

The truth of these lemmas follows easily from the relevant definitions.  Finally, we
should verify that the measure $\pi_* \nu$ is a suitable measure to work with.

\begin{lemma}
\label{aclem}
The measure $\pi_* \nu$ is an absolutely continuous invariant probability measure
for $f$, with respect to Lebesgue measure.
\end{lemma}

\emph{Proof:} Firstly, $\pi_* \nu$ is easily seen to be a mixing invariant
probability measure for $f$.
The measure $\nu$ is equivalent to $m$, so $\pi_* \nu$ is equivalent
to $\pi_* m$.  We show $\pi_* m \ll \mathrm{Leb}$.
Let $A \subset X$, and assume $\pi_* m(A)>0$.  Let
$\eta$ be the usual partition of $\Delta$; there exists $E \in \eta$ with
$m(\pi^{-1} A \cap E)>0$.  So there exists $E_0 \in
\mathcal{P}$ and $i<R(E_0)$ such that $m(F^{-i}
\pi^{-1} A \cap E_0)>0$.  This means there exists $A_0 \subset E_0$
such that $m(A_0)>0$ and $F^i A_0 \subset \pi^{-1} (A)$.

We know from the construction
that $\mathrm{Leb} (A_0) = m (A_0) > 0$.  Also, $f^i | A_0$ is injective, and since
the Jacobian of $f^R$ is non-singular, so is the Jacobian of $f$, and
$\mathrm{Leb}(f^i A_0)>0$.  We have
$f^i (A_0) = \pi F^i (A_0) \subset A$, so $\mathrm{Leb}(A)>0$.  This
proves that $\pi_* m \ll \mathrm{Leb}$.
\hfill $\square$

\section{An example with non-H\"older projection}
\label{baddist}

Finally, we briefly give an example of a system which exhibits a tower
with non-H\"older tower projection.

Given $0<a<b<1$, and $\alpha>1$, we define $f:[0,1] \circlearrowleft$ by
\begin{eqnarray*} f(x) = \left\{ \begin{array}{rl}
  1 - (1-b) (-\log a)^\alpha (- \log (a-x))^{-\alpha} & x \in [0,a] \\
 \frac{b}{b-a} (x-a) & x \in (a,b) \\
 b - \frac{b}{a} \exp \{ (1-b)^{\alpha^{-1}} (\log a) (1-x)^{-\alpha^{-1}} \}
   & x \in [b,1] \end{array} \right. .
\end{eqnarray*}
\begin{figure}
\begin{center}
\includegraphics[height=1.8in]{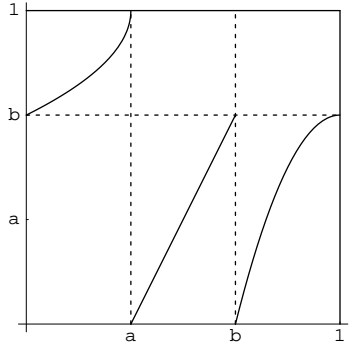}
\caption{\label{fig}A system admitting a tower via a non-H\"older semi-conjugacy.}
\end{center}
\end{figure}
(See figure \ref{fig}.)
Note that the map has unbounded derivative near $a$, and that $a$ maps onto
the
critical point at $1$.  It is easy to check that $f$ is monotone increasing on each
interval, and that $f$ has a Markov structure on the intervals $[0,a]$, $(a,b)$,
$[b,1]$, and could be represented by a subshift of finite type.
A straightforward, if tedious, calculation shows that for $x \in [0,a]$,
$y \in [b,1]$,
\begin{eqnarray*}
  f'(x) &=& \alpha (1-b) (-\log a)^\alpha (a-x)^{-1} (- \log (a-x))^{-(\alpha+1)} ;\\
  f'(y) &=& - \frac{b}{a} \alpha^{-1} (1-b)^{\alpha^{-1}} (\log a)
		(1-y)^{-(\alpha^{-1} + 1)} \\
       && \, \, \, \, \, \, \, \, \, \, \, \, 
		\, \, \, \,
		\exp \{ (1-b)^{\alpha^{-1}} (\log a)
			(1-y)^{-\alpha^{-1}} \} ;\\
  f'(f(x)) &=& \frac{b}{a} \alpha^{-1} (1-b)^{-1} (- \log a)^{-\alpha}
		(a-x)^{-1} (- \log (a-x))^{\alpha+1} ;\\
  (f^2)'(x) &=& \frac{b}{a} .
\end{eqnarray*}

Taking $\Delta_0 = [0,b]$, $\mathcal{P} = \{ [0,a],(a,b) \}$ with
$R([0,a])=2$, $R((a,b))=1$, it is clear that the conditions for $f$ to admit
a tower $F:\Delta \circlearrowleft$ are satisfied.
For $x,y \in [0,a]$,
we have that $|x-y| \approx (\frac{b}{a})^{-s(x,y)}$.
If we fix $k$ and consider $|f(x) - f(y)|$ for $x,y \in [0,a)$
with $s(x,y)=k$, then for $y$ close to $a$, we have
\begin{eqnarray*}
   |f(x) - f(y)| &\approx& (k \log (b/a) + C)^{-\alpha} \\
      &\approx& k^{-\alpha} \end{eqnarray*}
for some $C$.  This determines the regularity of the tower projection
$\pi$, which is in particular not H\"older continuous.  However, if we
take $\psi \in (R1,\gamma)$ for some $\gamma \in (0,1]$, then for
$x,y \in \Delta$,
\[ | \psi \circ \pi (x) - \psi \circ \pi (y)| \leq C s(x,y)^{-\alpha
   \gamma} ,\]
so $\psi \circ \pi \in (V4,\alpha \gamma)$.  This means that for
observables in $(R1,\gamma)$ for sufficiently large $\gamma$ (dependent
on $\mathcal{R}_\infty (\psi)$), we get polynomial decay of correlations.
Also note that if we take $\psi$ with 
$\mathcal{R}_\varepsilon (\psi) = \mathcal{O}
  (\theta^{\varepsilon^{-\alpha^{-1}}})$,
for some $\theta \in (0,1)$, then
\[ | \psi \circ \pi (x) - \psi \circ \pi (y)| \leq C \theta^{(C s(x,y))} ,\]
and $\psi \circ \pi \in (V1)$.  So we get exponential mixing for some
non-trivial subset of $(R1)$.  We may find observable classes corresponding
to $(V2)$, $(V3)$ in similar fashion.

\bibliographystyle{plain}

\begin{thebibliography}{99}

\bibitem[ALP]{ALP1} J. F. Alves, S. Luzzatto, V. Pinheiro, \emph{Markov
structures and decay of correlations for non-uniformly expanding dynamical
systems}, preprint (2002).
\bibitem[AV]{AV} J. F. Alves, M. Viana, \emph{Statistical stability for
robust classes of maps with non-uniform expansion}, Ergod. Th. Dynam.
Sys. {\bf 22}, no. 1, 1-32 (2002).
\bibitem[Ba]{Ba} V. Baladi, \emph{Positive transfer operators and
decay of correlations}, Advanced Series in Nonlinear Dynamics, {\bf 16},
World Scientific Publishing Co., Inc., River Edge, NJ (2000).
\bibitem[Bo]{B} R. Bowen, \emph{Equilibrium states and the ergodic theory
of Anosov diffeomorphisms}, Springer Lecture Notes in Math., {\bf 470}
(1975).
\bibitem[BFG]{BFG} X. Bressaud, R. Fern\'andez, A. Galves, \emph{Decay of
correlations for non-H\"olderian dynamics. A coupling approach}, Electron.
J. Probab. {\bf 4}, no. 3, 19 pp. (electronic) (1999).
\bibitem[BLS]{BLvS} H. Bruin, S. Luzzatto, S. van Strien, \emph{Decay of
correlations in one-dimensional dynamics}, to appear,
Annales de l'ENS.
\bibitem[BM1]{BM} J. Buzzi, V. Maume-Deschamps, {\em Decay of correlations
on towers with non-H\"older Jacobian and non-exponential return time},
preprint (2001).
\bibitem[BM2]{BM2} J. Buzzi, V. Maume-Deschamps, {\em Decay of correlations
for piecewise invertible maps in higher dimensions}, Israel. J. Math,
{\bf 131}, 203-220 (2002).
\bibitem[BST]{BST} J. Buzzi, O. Sester, M. Tsujii, \emph{Weakly expanding
   skew-products of quadratic maps}, to appear, Ergod. Th. Dynam. Sys.
\bibitem[FL]{FL} A. Fisher, A. Lopes, {\em Exact bounds for the polynomial
decay of correlation, $1 / f$ noise and the CLT for the equilibrium state of
a non-H\"older potential}, Nonlinearity, {\bf 14}, 1071-1104 (2001).
\bibitem[Hol]{H1} M. Holland, {\em Physical measures for chaotic dynamical
systems and decay of correlations}, Ph.D. thesis (Warwick, 2001).
\bibitem[Hu]{Hu} H. Hu, \emph{Decay of correlations for piecewise smooth
   maps with indifferent fixed points}, preprint.
\bibitem[I]{I} S. Isola, \emph{On systems with finite ergodic degree},
   preprint.
\bibitem[KMS]{KMS} A. Kondah, V. Maume, B. Schmitt, \emph{Vitesse de
convergence vers l'\'etat d'\'equilibre pour des dynamiques markoviennes
non h\"old\'eriennes}, Ann. Inst. H. Poincar\'e Probab. Statist., {\bf 33},
no. 6, 675-695 (1997).
\bibitem[L]{L2} C. Liverani, \emph{Central limit theorem for deterministic
systems}, International conference on dynamical systems, Montevideo 1995,
Eds. F. Ledrappier, J. Lewowicz, S. Newhouse, Pitman research notes in Math.
{\bf 362}, 56-75 (1996).
\bibitem[LM]{LM} M. Lyubich, J. Milnor, \emph{The Fibonacci unimodal map},
J. Amer. Math. Soc. {\bf 6}, 425-457 (1993).
\bibitem[M]{M} V. Maume-Deschamps, \emph{Projective metrics and mixing
properties on towers}, Trans. Amer. Math. Soc., {\bf 353}, 3371-3389 (2001).
\bibitem[P]{P} M. Pollicott, \emph{Rates of mixing for potentials of
   summable variation}, Trans. Amer. Math. Soc, {\bf 352}, 843-853, (2000).
\bibitem[PY]{PY} M. Pollicott, M. Yuri, \emph{Statistical Properties of Maps
   with Indifferent Periodic Points}, Commun. Math. Phys. {\bf 217},
   503-520 (2001).
\bibitem[R]{R} D. Ruelle, A measure associated with axiom-A attractors,
Amer. J. Math. {\bf 98}, 619-654 (1976).
\bibitem[Si]{Si} Ja. G. Sina\u\i, Gibbs measures in ergodic theory, Uspehi
   Mat. Nank {\bf 27}, 21-64 (1972) (Russian).
\bibitem[V1]{V} M. Viana, \emph{Multidimensional nonhyperbolic attractors},
   Inst. Hautes \'Etudes Sci. Publ. Math. {\bf 85}, 63-96 (1997).
\bibitem[V2]{V2} M. Viana, \emph{Stochastic dynamics of deterministic
systems}, Brazillian Math. Colloquium, IMPA (1997).
\bibitem[Yo]{Y2} L.-S. Young,
{\em Recurrence times and rates of mixing}, Israel. J. of Math. {\bf 110},
153-188 (1999).

\end{thebibliography}

\end{document}